\documentclass[12pt]{article}

\usepackage{a4wide}
\usepackage{amssymb,amsmath,array}

\allowdisplaybreaks[1]

\newcommand{\titre}{Primitive ideals and automorphisms\\ of quantum matrices}

\newenvironment{proof}{\begin{trivlist}\item[]{\it
Proof.}}{\hfill$\square$\end{trivlist}}

\newtheorem{theorem}{Theorem}[section]
\newtheorem{corollary}[theorem]{Corollary}

\newtheorem{lemma}[theorem]{Lemma}
\newtheorem{proposition}[theorem]{Proposition}

\newtheorem{claim}[theorem]{Claim}

\newcommand{\gc}{ [ \hspace{-0.65mm} [}
\newcommand{\dc}{]  \hspace{-0.65mm} ]}
\newcommand{\glm}{\mathrm{GL}_{m}(\mathbb{C}) }

\newcommand{\ia}{i,\alpha}

\newcommand{\fract}{\mathrm{Frac}}
\newcommand{\rank}{\mathrm{rk}}

\newcommand{\hc}{\mathcal{H}}

\newcommand{\diag}{{\rm diag}}

\newcommand{\spec}{{\rm Spec}}
\newcommand{\aut}{{\rm Aut}}
\newcommand{\prim}{{\rm Prim}}

\newcommand{\ideal}[1]{\left\langle {#1} \right\rangle}

\def\mc{{\mathbb{C}}}
\def\N{{\mathbb{N}}}

\def\detq{{\rm det}_q}

\def\ch{{\mathcal H}}
\def\co{{\mathcal O}}

\def\oq{{\cal O}_q}
\def\oqmn{\co_q(M_n)}

\def\oqmtwo{\co_q(M_2)}
\def\oqmm13{\co_q(M_{1,3})}
\def\oqm23{\co_q(M_{2,3})}
\def\oqmmn{\co_q(M_{m,n})}
\def\oqmnm{\co_q(M_{n,m})}

\def\qdot{q^{\bullet}}

\input{xy}
\xyoption{all}

\begin{document}

\title{\titre}
\author{S Launois and T H Lenagan
\thanks{This research was supported by a Marie Curie Intra-European
  Fellowship within the $6^{\mbox{th}}$ European Community Framework
  Programme and by Leverhulme Research Interchange
Grant F/00158/X }\;
}
\date{}

\maketitle


\begin{abstract}
Let $q $ be a nonzero complex number that is not a root of unity. We give a
criterion for $\ideal{0}$ to be a primitive ideal of the algebra $\oqmmn$ of
quantum matrices. Next, we describe all height one primes of $\oqmmn$; these
two problems are actually interlinked since it turns out that $\ideal{0}$ is a
primitive ideal of $\oqmmn$ whenever $\oqmmn$ has only finitely many height
one primes. Finally, we compute the automorphism group of $\oqmmn$ in the case
where $m \neq n$. In order to do this, we first study the action of this group
on the prime spectrum of $\oqmmn$. Then, by using the preferred basis of
$\oqmmn$ and PBW bases, we prove that the automorphism group of $\oqmmn$ is
isomorphic to the torus $(\mathbb{C}^*)^{m+n-1}$ when $m \neq n$ and $(m,n)
\neq (1,3),(3,1)$.

\end{abstract}



\section*{Introduction}\label{sec:intro}

The automorphism group of a polynomial algebra $\mathbb{C}[X_1,\dots,X_n]$,
for $n \geq 2$, is in general far from being understood. 
In the case $n=2$  
the group was described explicitly by Jung, \cite{jung}. 
However, it
is only recently that Shestakov and Umirbaev, \cite{umirbaev}, have proved
that the well-known Nagata automorphism of $\mathbb{C}[X_1,X_2,X_3]$ is wild;
that is, the Nagata automorphism cannot be written as a product of elementary
automorphisms.

In this paper, we are interested in the quantum case. More precisely, when $q
$ is a nonzero complex number that is not a root of unity, we study the
automorphism group of the algebra $\oqmmn$ of $m \times n$ quantum matrices
which in turn is a non-commutative deformation of a polynomial ring in $m
\times n$ indeterminates. We denote by $Y_{\ia}$, $(i,\alpha) \in \gc 1,m \dc
\times \gc 1,n \dc$, the canonical generators of $\oqmmn$. It is well-known
that the group $\ch:=(\mathbb{C}^*)^{m+n}$ acts on $\oqmmn$ by
$\mathbb{C}$-automorphisms via: $$(a_1,\dots,a_m; b_1,\dots,b_n).Y_{\ia}=a_i
b_{\alpha} Y_{\ia} \quad ((\ia) \in \gc 1,m \dc \times \gc 1,n \dc).$$ (Note
that this is not a faithful action; for example, $(a_,\dots,a; 1,\dots,1)$ and
$(1,\dots, 1;a,\dots,a)$ have the same action on $\oqmmn$, multiplying each
$Y_{\ia}$ by $a$. This explains why the automorphism group is an
$(m+n-1)$-torus rather than an $(m+n)$-torus.)

It was observed by Alev and Chamarie, \cite{alevchamarie}, that quantization
implies rigidity and so puts limits on the automorphism group of quantized
algebras. This explains why it has been possible to compute the automorphism
group of several quantum algebras at least in the generic case. For instance,
Alev and Chamarie, \cite{alevchamarie}, have described the automorphism group
of the algebra of $2 \times 2$ quantum matrices,
Alev and Dumas, \cite{alevdumasrigidite}, the automorphism group of the
positive part of the quantized enveloping algebra of a complex simple Lie
algebra of type $A_2$, the first author, \cite{lau}, the automorphism group of
the positive part of the quantized enveloping algebra of a complex simple Lie
algebra of type $B_2$, ... .

One method to study the automorphism group of an algebra 
is to use the invariance of the set
of height one primes under the action of the automorphism group.
This method was used successfully by Rigal, \cite{rigal}, to compute the
automorphism group of quantized Weyl algebras and next by Gomez-Torrecillas
and El Kaoutit, \cite{gomes}, to calculate the automorphism group of the
coordinate ring of quantum symplectic spaces. 
We also use this method 
in the present paper. However, in the cases of quantized Weyl algebras and
coordinate ring of quantum symplectic spaces, the number of height one primes
is finite (because of the choice of parameters) and so the restrictions
on the automorphisms are very strong. In the case of the algebra of quantum
matrices, the set of height one primes is in general not finite; so the
situation is substantially more complicated.

We start by considering when the algebra $\oqmmn$ has only finitely many
height one primes. It turns out that, because of the stratification theorem of
Goodearl and Letzter, see \cite{bg} for example, this situation arises exactly
when $\ideal{0}$ is a primitive ideal of $\oqmmn$. Thus, the first section of
this paper is devoted to this question: we establish a criterion for the ideal
$\ideal{0}$ to be primitive in $\oqmmn$. More precisely, we prove that
$\ideal{0}$ is a primitive ideal of $\oqmmn$ if and only if $v_2(m) \neq
v_2(n)$, where $v_2(k)$ denotes the $2$-adic valuation of a positive integer
$k$. This shows, for example, that $\ideal{0}$ is a primitive ideal of
$\oqm23$. This criterion together with the stratification theorem of Goodearl
and Letzter shows that, if $v_2(m) \neq v_2(n)$, then $\oqmmn$ has only
finitely many height one primes and it turns out that they are all
$\ch$-invariant. On the other hand, if $v_2(m) = v_2(n)$, then $\oqmmn$ has
infinitely many height one primes. A finite number (those that are
$\ch$-invariant) are already known from 
previous work of the authors and Rigal,
\cite{llr}. In the second part of this paper, we provide an explicit generator
for every height one prime of $\oqmmn$. (Note that, since $\oqmmn$ is a
Noetherian (non-commutative) UFD, every height one prime of $\oqmmn$ is
generated by a normal element, \cite{llr}.) Hence, in this second part, we
describe all normal elements of $\oqmmn$ that generate a prime ideal.

Finally, in the third section, we investigate the automorphism group of
$\oqmmn$ when $m \neq n$. Using the description of the height one primes of
$\oqmmn$ together with graded arguments, we show that all height one primes
that are $\ch$-invariant, except possibly one, are invariant under every
automorphism of $\oqmmn$. Next, by using the preferred basis of $\oqmmn$ and
certain PBW bases, we are able to prove that the automorphism group of
$\oqmmn$ is isomorphic to the torus $(\mathbb{C}^*)^{m+n-1}$ when $m \neq n$
and $(m,n) \neq (1,3),(3,1)$. This latter restriction occurs because of the
existence of a rogue automorphism for quantum $3$-space (which may be viewed
as $1\times 3$ quantum matrices), see \cite{alevchamarie}. This is the only
exception, and our analysis recovers the Alev-Chamarie result for this
exceptional case.

In the case where $m=n$, the algebra $O_q(M_n)$ has a homogeneous central
element of degree $n$, the quantum determinant. The existence of this element
considerably complicates the computation of the automorphism group, since the
graded arguments that we use in the non-squared case do not put strong limits
on the automorphism group and its action on the height one primes in the
square case. In addition, transposition provides an automorphism that also
clouds the analysis. The conjecture is that the automorphim group of $\oqmn$
is generated by a torus and the transposition automorphism. This has been
verified by Alev and Chamarie, \cite{alevchamarie}, in the $2\times 2$ case.
Our present methods can recover the Alev-Chamarie result, but, as yet, we are
unable to deal with the general case, although we do have partial results. We
intend to return to this question in a subsequent paper.


\section{A criterion for quantum matrices to be primitive.}

In this section, we use the $\ch$-stratification theory of Goodearl and
Letzter together with the deleting derivations theory of Cauchon in order to
characterize the integers $m$ and $n$ such that $\ideal{0}$ is a primitive
ideal in $\oqmmn$.


\subsection{The $\ch$-stratification of the prime spectrum of $\oqmmn$.}
$ $

Throughout this paper, we use the following conventions. \\$ $ 
\\$\bullet$ If
$I$ is a finite set, $|I|$ denotes its cardinality.  
\\$\bullet$  $\gc a,b \dc := \{ i\in{\mathbb N} \mid a\leq i\leq b\}$. 
\\$\bullet$ $\mathbb{C}$
denotes the field of complex numbers and we set
$\mathbb{C}^*:=\mathbb{C}\setminus \{0\}$.  
\\$\bullet$ \textbf{$q\in
\mathbb{C}^*$ is not a root of unity}. 
\\$\bullet$ $m,n$ denote positive
integers. 
\\$\bullet$ $R=\oqmmn$ is the quantization of the ring of
regular functions on $m \times n$ matrices with entries in $\mathbb{C}$; it is
the $\mathbb{C}$-algebra generated by the $m \times n $ indeterminates
$Y_{\ia}$, $1 \leq i \leq m$ and $ 1 \leq \alpha \leq n$, subject to the
following relations:\\ \[
\begin{array}{ll}
Y_{i, \beta}Y_{i, \alpha}=q^{-1} Y_{i, \alpha}Y_{i ,\beta},
& (\alpha < \beta); \\
Y_{j, \alpha}Y_{i, \alpha}=q^{-1}Y_{i, \alpha}Y_{j, \alpha},
& (i<j); \\
Y_{j,\beta}Y_{i, \alpha}=Y_{i, \alpha}Y_{j,\beta},
& (i <j,  \alpha > \beta); \\
Y_{j,\beta}Y_{i, \alpha}=Y_{i, \alpha} Y_{j,\beta}-(q-q^{-1})Y_{i,\beta}Y_{j,\alpha},
& (i<j,  \alpha <\beta). 
\end{array}
\]

It is well-known that $R$ can be presented as an iterated Ore extension over
$\mathbb{C}$, with the generators $Y_{\ia}$ adjoined in lexicographic order.
Thus the ring $R$ is a Noetherian domain; we denote by $F$ its skew-field of
fractions. Moreover, since $q$ is not a root of unity, it follows from
\cite[Theorem 3.2]{gletpams} that all prime ideals of $R$ are completely
prime. We denote by $\spec(R)$ the set of (completely) prime ideals of $R$.
\\$\bullet$ It is well-known that the algebras $\oqmmn$ and $\oqmnm$ are
isomorphic. Hence, we assume that $m \leq n$. \\$\bullet$ It is easy to check
that the group $\ch:=\left( \mathbb{C}^* \right)^{m+n}$ acts on $R$ by
$\mathbb{C}$-algebra automorphisms via:
$$(a_1,\dots,a_m,b_1,\dots,b_n).Y_{\ia} = a_i b_\alpha Y_{\ia} \quad {\rm
for~all} \quad \: (\ia)\in \gc 1,m \dc \times \gc 1,n \dc.$$ An $\ch
$-eigenvector $x$ of $R$ is a nonzero element $x \in R$ such that $h.x \in
\mathbb{C}^*x$ for each $h \in \ch$. An ideal $I$ of $R$ is said to be
$\ch$-invariant if $h.I =I$ for all $h\in \ch$. We denote by $\ch$-$\spec(R)$
the set of $\ch$-invariant prime ideals of $R$. Since $q$ is not a root of
unity, it follows from \cite[5.7]{glet} that $\ch$-$\spec(R)$ is a finite set.
Note that $\ideal{0}$ is an $\ch$-invariant prime
ideal of $R$,  since $R$ is a domain. \\$ $

The action of $\ch$ on $R$ allows  us to use 
the $\ch$-stratification theory of
Goodearl and Letzter, see \cite[II.2]{bg}, to constuct a partition of
$\spec(R)$ as follows. If $J$ is an $\ch$-invariant prime ideal of $R$, we
denote by $\spec_J(R)$ the $\ch$-stratum of $\spec(R)$ associated to $J$.
Recall that $\spec_J(R):=\{ P \in \spec(R) \mid \bigcap_{h \in \ch} h.P=J \}$.
Then the $\ch$-strata $\spec_J(R)$, 
with $J \in \ch$-$\spec(R)$, form a partition
of $\spec(R)$, see \cite{bg}:
\begin{eqnarray}
\label{partition}
\spec(R) & = & \bigsqcup_{J \in \ch \mbox{-}\spec(R)} \spec_J(R).
\end{eqnarray}
Naturally, this partition induces a partition of the set $\prim(R)$ of all
(left) primitive ideals of $R$ as follows. For all $J \in \ch$-$\spec(R)$, we
set $\prim_J(R):=\spec_J(R) \cap \prim(R)$. Then it is obvious that the
$\ch$-strata $\prim_J(R)$ ($J \in \ch$-$\spec(R)$) form a partition of
$\prim(R)$: 
$$
\prim(R) = \bigsqcup_{J \in \hc \mbox{-}\spec(R)} \prim_J(R).
$$
One of the reasons that makes the $\ch$-stratification interesting is that it
provides a powerful tool for recognizing primitive ideals. Indeed, since
$\mathbb{C}$ is uncountable and since the Noetherian domain $R$ is generated
as an algebra by a finite number of elements, it follows from
\cite[Proposition II.7.16]{bg} that the algebra $R$ satisfies the
Nullstellensatz over $\mathbb{C}$, see \cite[II.7.14]{bg}. Further the set of
$\ch$-invariant prime ideals of $R$ is finite. Thus \cite[Theorem II.8.4]{bg}
implies that $\prim_J(R)$ $(J \in \hc$-$\spec(R)$) coincides with the set of
those primes in $\spec_J(R)$ that are maximal in $\spec_J(R)$. Hence, we
deduce the following criterion for $R$ to be a primitive ring.

\begin{proposition} \label{primitivity1}
$R$ is a primitive ring if and only if the $\ideal{0}$-stratum of $\spec(R)$
is reduced to $\ideal{0}$; that is, $\ideal{0}$ is the only prime ideal in 
 $\spec_{\ideal{0}}(R)$.
\end{proposition}



\subsection{Dimension of the $\ideal{0}$-stratum of $\spec(R)$.}
\label{subsectionStrate0}

In this section, we show that the $\ideal{0}$-stratum $ \spec_{\ideal{0}}
(\oqmmn)$ is Zariski-homeomorphic to the prime spectrum of a commutative
Laurent polynomial ring in $\dim (\ker(B))$ indeterminates, where $B$ is a $mn
\times mn$ matrix with entries in $\mathbb{C}$. First, we describe
the matrix $B$ explicitly and then we will compute the dimension of the kernel of $B$. \\

First recall, see \cite{c2}, that the theory of deleting derivations can be
applied to the iterated Ore extension $R=\mathbb{C}[Y_{1,1}]\dots
[Y_{m,n};\sigma_{m,n},\delta_{m,n}]$ (where the indices are increasing for the
lexicographic order $\leq$). The corresponding deleting derivations algorithm
is called the standard deleting derivations algorithm. Before recalling its
construction, we need to introduce some notation.

\begin{itemize}
\item We denote by $\leq_s$ the lexicographic ordering on $\mathbb{N}^2$. We
often call it the standard ordering on $\mathbb{N}^2$. Recall that $(\ia)
\leq_s (j,\beta)$ if and only if $ [(i < j) \mbox{ or } (i=j \mbox{ and }
\alpha \leq \beta )]$. 
\item We set $E_s=\left(\gc 1,m \dc \times \gc 1,n \dc
\cup \{(m,n+1)\} \right) \setminus \{(1,1)\}$. \item Let $(j,\beta) \in E_s$.
If $(j,\beta) \neq (m,n+1)$, then $(j,\beta)^{+}$ denotes the least element
(relative to $\leq_s$) of the set $\left\{ (\ia) \in E_s \mbox{ $\mid$
}(j,\beta) <_s (\ia) \right\}$.
\end{itemize}

As described in \cite{c2}, 
the standard deleting derivations algorithm  
constructs, for each $r \in E_s$, a family $(Y_{\ia}^{(r)})_{(\ia) \in
\gc 1,m \dc \times \gc 1,n \dc}$ of elements of $F:=\fract(R)$, defined as
follows. \\$ $
\begin{enumerate}
\item \underline{If $r=(m,n+1)$}, then
$Y_{\ia}^{(m,n+1)}=Y_{\ia}$ for all $(\ia) \in \gc 1,m \dc \times \gc 1,n
\dc$.\\$ $
\item \underline{Assume that $r=(j,\beta) <_s (m,n+1)$}
and that the $Y_{\ia}^{(r^{+})}$ ($(\ia) \in \gc 1,m \dc \times \gc 1,n\dc$)
are already constructed.
Then, it follows from \cite[Th\'eor\`eme 3.2.1]{cauchoneff} that
$Y_{j,\beta}^{(r^+)} \neq 0$ and,
for all $(\ia) \in \gc 1,m \dc \times \gc 1,n\dc$, we have:
$$Y_{\ia}^{(r)}=\left\{ \begin{array}{ll}
Y_{\ia}^{(r^{+})}-Y_{i,\beta}^{(r^{+})}
\left(Y_{j,\beta}^{(r^{+})}\right)^{-1}
Y_{j,\alpha}^{(r^{+})}
& \mbox{ if } i<j \mbox{ and } \alpha < \beta \\
Y_{\ia}^{(r^{+})} & \mbox{ otherwise.}
\end{array} \right.$$
\end{enumerate}

As in \cite{cauchoneff}, we denote by $\overline{R}$ the subalgebra of
$\fract(R)$ generated by the indeterminates obtained at the end of this
algorithm, that is, we denote by $\overline{R}$ the subalgebra of $\fract(R)$
generated by the $T_{\ia}:=Y_{\ia}^{(1,2)}$ for each 
$(\ia) \in \gc 1,m \dc \times \gc
1,n \dc$. \\

Let $N \in \mathbb{N}^*$ and let $\Lambda=(\Lambda_{i,j})$ 
be a multiplicatively antisymmetric $N\times N$ 
matrix over $\mathbb{C}^*$; 
that is, $\Lambda_{i,i}=1$ and $\Lambda_{j,i}=\Lambda_{i,j}^{-1}$ 
for all $i,j \in \gc 1,N \dc$. 
We denote by $\mathbb{C}_{\Lambda}[T_1,\dots,T_N]$ 
the corresponding quantum affine space; that is, the $\mathbb{C}$-algebra 
generated by the $N$ indeterminates $T_1,\dots,T_N$ subject to the 
relations $T_i T_j =\Lambda_{i,j} T_j T_i $ for all $i,j \in \gc 1,N \dc$. 
In \cite[Section 2.2]{c2}, Cauchon has shown that $\overline{R}$ can be 
viewed as the quantum affine space generated by the indeterminates 
$T_{\ia}$  for $(\ia) \in \gc 1,n \dc^2$, subject to the following 
relations. 
\[
\begin{array}{ll}
T_{i, \beta}T_{i, \alpha}=q^{-1}T_{i, \alpha}T_{i ,\beta},
& (\alpha < \beta); \\
T_{j, \alpha}T_{i, \alpha}=q^{-1}T_{i, \alpha}T_{j, \alpha},
& (i<j); \\
T_{j,\beta}T_{i, \alpha}=T_{i, \alpha}T_{j,\beta},
& (i <j,  \alpha > \beta); \\
T_{j,\beta}T_{i, \alpha}=T_{i, \alpha} T_{j,\beta},
& (i<j,  \alpha <\beta). 
\end{array}
\]

Hence $\overline{R}=\mathbb{C}_{\Lambda}[T_{1,1},T_{1,2},\dots,T_{m,n}]$,
where $\Lambda$ denotes the $mn \times mn$ matrix defined as follows. We set
$$A:=\left(
\begin{array}{ccccc}
 0 & 1 & 1 & \dots & 1 \\
-1 & 0 & 1 & \dots  & 1 \\
\vdots & \ddots &\ddots&\ddots &\vdots \\
-1 & \dots & -1 & 0 & 1  \\
 -1& \dots& \dots & -1 & 0 \\  
\end{array} 
\right) 
\in \mathcal{M}_{m}(\mathbb{C}), $$
and 
 $$B:=\left( 
\begin{array}{ccccc}
 A & I_m & I_m & \dots & I_m \\
-I_m & A & I_m & \dots  & I_m \\
\vdots & \ddots &\ddots&\ddots &\vdots \\
-I_m & \dots & -I_m & A & I_m  \\
 -I_m& \dots& \dots & -I_m & A \\  
\end{array} 
\right)
\in \mathcal{M}_{mn}(\mathbb{C}), $$
where $I_m$ denotes the identity matrix of $\mathcal{M}_m$. Then $\Lambda $ is
the $mn \times mn$ matrix whose entries are defined by $\Lambda_{k,l}
=q^{b_{k,l}}$ for all $k,l \in \gc 1, mn \dc$.

It follows from
\cite[Th\'eor\`emes 5.1.1, 5.5.1, 5.5.2]{cauchoneff} that $\spec_{\ideal{0}}
(R)$ is Zariski-homeomorphic to
the prime spectrum of the quantum torus $P(\Lambda):=\overline{R}
\Sigma^{-1}$, where $\Sigma$ denotes the multiplicative system of
$\overline{R}$ generated by the normal elements $T_{\ia}$ with 
$(\ia) \in \gc 1,m
\dc \times \gc 1,n \dc$.
Next, 
$\spec(P(\Lambda))$ is Zariski-homeomorphic via extension and contraction to
the prime spectrum of the centre $Z(P(\Lambda))$ of $P(\Lambda)$, by
\cite[Corollary 1.5]{glettore}.
Further,
$Z(P(\Lambda))$ turns out to be a Laurent polynomial ring. To make 
this result precise, we need to introduce the following notation.

If $\underline{s}=(s_{1,1},s_{1,2},\dots,s_{m,n}) \in \mathbb{Z}^{mn}$, then
we set $T^{\underline{s}}:= T_{1,1}^{s_{1,1}}T_{1,2}^{s_{1,2}} \dots
T_{m,n}^{s_{m,n}} \in P( \Lambda)$.

As in \cite{glettore}, we denote by $\sigma : \mathbb{Z}^{mn} \times
\mathbb{Z}^{mn} \rightarrow \mathbb{C}^*$ the antisymmetric bicharacter
defined by 
$$\sigma(\underline{s},\underline{t}):=\prod_{k,l=1}^{mn}
\Lambda_{k,l}^{s_kt_l} \: \quad {\rm for~all}\quad  
\underline{s},\underline{t} \in
\mathbb{Z}^{mn}.$$ 
Then it follows from \cite[1.3]{glettore} that the centre
$Z(P(\Lambda))$ of $P(\Lambda)$ is a Laurent polynomial ring in the variables
$(T^{\underline{b_1}})^{\pm 1},\dots,(T^{\underline{b_r}})^{\pm 1}$, where
$(\underline{b_1},\dots,\underline{b_r})$ is any basis of $S:=\{\underline{s}
\in \mathbb{Z}^{mn} \mid \sigma(\underline{s},-)\equiv 1\}$. Since $q$ is not
a root of unity, easy computations show that $\underline{s} \in S$ if and only
if $B^t\underline{s}^t=0$. Hence the centre $Z(P(\Lambda))$ of $P(\Lambda)$
is a Laurent polynomial ring in $\ker(B^t)$ indeterminates; so we have just
proved the following statement.

\begin{proposition} \label{dimension0stratum}
$\spec_{\ideal{0}} (\oqmmn)$ is Zariski-homeomorphic to the prime spectrum of
a commutative Laurent polynomial ring in $\dim (\ker(B^t))=\dim (\ker(B))$
indeterminates. 
\end{proposition}

We deduce from Propositions \ref{primitivity1} and \ref{dimension0stratum} the
following criterion for $R$ to be a primitive ring.

\begin{corollary} \label{primitivity2}
$R$ is a primitive ring if and only if $B$ is invertible.
\end{corollary}

We now compute the dimension of the kernel of $B$. First, straightforward computations show that 
$$
\dim(\ker(B)) = m- \rank \left[ (A+I)^n +
(A-I)^n \right].
$$

If $i$ is a positive integer greater than or equal to $2$, we denote by
$v_2(i)$ the 2-adic valuation of $i$; that is, $m$ is the largest integer such
that $2^m|i$.

One can easily check the following result.

\begin{proposition} \label{theoremdimkerB}
$$
\dim(\ker(B))=\left\{ \begin{array}{cc}
0 & \mbox{ if $v_2(m) \neq v_2(n)$} \\
m\wedge n & \mbox{ otherwise} \\
\end{array} \right.  
$$
where $m \wedge n$ denotes the greatest common divisor of $m$ and $n$.
\end{proposition}

\begin{proof} 
Set 
$e_k := \exp\left(i\left( \frac{(2k+1)\pi }{m} \right) \right)$ for all $k
\in \gc 0,m-1 \dc$; that is, the $e_k$ are the $m$th root of $-1$. It is easy to show that there exists $U \in \glm$ such that $U^{-1} A U =
\diag(\frac{e_0+1}{e_0-1},\dots,\frac{e_{m-1}+1}{e_{m-1}-1})$. 
Hence we have
\begin{eqnarray*}
\lefteqn{U^{-1}
\left[(A+I)^n + (A-I)^n \right] U = }\\[2ex]
&&\diag\left( \left(\frac{2e_0}{e_0-1}
\right)^n + \left(\frac{2}{e_0-1} \right)^n,\dots ,
\left(\frac{2e_{m-1}}{e_{m-1}-1} \right)^n + \left(\frac{2}{e_{m-1}-1}
\right)^n \right).
\end{eqnarray*} 
Since $U \in \glm$, 

\begin{eqnarray*}
\rank \left( (A+I)^n + (A-I)^n \right) & = & \rank \left( U^{-1} \left[
(A+I)^n + (A-I)^n \right] U \right) \\ & = & \begin{array}{|c|} \{k \in \gc
0,m-1 \dc \mid \left(\frac{2e_k}{e_k-1} \right)^n + \left(\frac{2}{e_k-1}
\right)^n \neq 0 \} 
\end{array},
\end{eqnarray*}
so that
$$
\dim(\ker (B))= \begin{array}{|c|} \left\{k \in \gc 0,m-1 \dc \mid
\left(\frac{2e_k}{e_k-1} \right)^n + \left(\frac{2}{e_k-1} \right)^n = 0
\right\}
\end{array} .
$$

Now, easy calculations show that $\left(\frac{2e_k}{e_k-1} \right)^n +
\left(\frac{2}{e_k-1} \right)^n = 0 $ if and only if $e_k^n=-1$; that is, if
and only if $\exp\left(i\left( \frac{(2k+1)n\pi }{m} \right) \right)=-1$, that
is, if and only if $m| (2k+1)n$ and $\frac{(2k+1)n }{m}$ is an odd integer. 
Hence, 
$$
\dim(\ker (B))= \begin{array}{|c|} \{k \in \gc 0,m-1 \dc \mid m| (2k+1)n
\mbox{ and } \frac{(2k+1)n }{m} {\rm is~odd} \} \end{array}.
$$

We now distinguish between two cases.
\\
$\bullet$ First, assume that $v_2(m) \neq v_2(n)$. Suppose that 
there exists $k \in \gc 0, m-1
\dc$ such that
$m| (2k+1)n$. Then $v_2(m) \leq v_2(n)$. 
This implies that $v_2(m) < v_2(n)$, since 
$v_2(m) \neq
v_2(n)$; and so $\frac{(2k+1)n }{m}$ is even. 
Hence there does not exist $k \in \gc 0,m-1 \dc$ such that $m|
(2k+1)n$ and $\frac{(2k+1)n }{m}$ is odd . This proves that, in this case, 
$$
\dim(\ker (B))= \begin{array}{|c|} \{k \in \gc 0,m-1 \dc \mid m| (2k+1)n
\mbox{ and } \frac{(2k+1)n }{m} \mbox{ is odd } \} \end{array} =0,
$$ 
as
desired.
\\$\bullet$ Next, assume that $v_2(m) = v_2(n)=\alpha$. Write
$n=2^{\alpha}v$ and $m = 2^{\alpha}u$ with
$u,v$ odd. It follows that 
\begin{eqnarray*}
\dim(\ker (B)) & = & 
\begin{array}{|c|} \{k \in \gc 0,m-1 \dc \mid u| (2k+1)v
\mbox{ and } \frac{(2k+1)v }{u} \mbox{ is odd } \} 
\end{array} \\ & = &
\begin{array}
{|c|} \{k \in \gc 0,m-1 \dc \mid u| (2k+1)v \} 
\end{array} .
\end{eqnarray*}
Set $d'=u\wedge v$ and $u=d'u'$. Note that
$d'$ 
and $u'$ are odd, since $u$ and $v$ are odd. 
This implies that $\dim(\ker (B))= \begin{array}{|c|}
\{k \in \gc 0,m-1 \dc \mid u' | (2k+1) \} \end{array}$, the number of odd
multiples of $u'$ less that $2m-1=2^{\alpha +1}d'u'-1$. Hence $\dim(\ker
(B))=2^{\alpha}d'=m \wedge n$, as required.
\end{proof}

Proposition \ref{theoremdimkerB} together with Proposition
\ref{dimension0stratum} and Corollary \ref{primitivity2} yields the following
results.

\begin{theorem} \label{theorem0stratum1}
$\spec_{\ideal{0}} (\oqmmn)$ is Zariski-homeomorphic to the prime spectrum of
a commutative Laurent polynomial ring in $\alpha_{m,n}$ indeterminates, where
$$
\alpha_{m,n}=\left\{ \begin{array}{cc} 0 & \mbox{ if $v_2(m) \neq v_2(n)$}
\\ m\wedge n & \mbox{ otherwise.} \\
\end{array} \right.
$$ 
\end{theorem}

\begin{theorem} \label{theorem0stratum2}
$\oqmmn$ is a primitive ring if and only if $v_2(m) \neq v_2(n)$.
\end{theorem}

For example, $\oqmtwo$ is not primitive, while $\oqm23$ is primitive.


\section{Height one primes in quantum matrices.}

In this section, we investigate height one primes of $R=\oqmmn$. 
Every height one prime ideal 
of $R$ is generated
by a normal element, since $R$ is
a Noetherian UFD, \cite{llr}. In this section, we describe explicitly the 
normal
elements that generate the height one 
prime ideals. Some of them are already known. Indeed,
the height one primes that are $\ch$-invariant have been described in
\cite{llr}. Hence, we mainly focus our attention on the other height one
primes: these belong to the $\ideal{0}$-stratum of $\spec(\oqmmn)$.

\subsection{Height one primes of $\oqmmn$ that are $\ch$-invariant.}

The algebra $\oqmn$ has a special element, $\detq$, the {\em quantum
determinant}, defined by 
\[ 
\detq := \sum_{\sigma}\,
(-q)^{l(\sigma)}Y_{1,\sigma(1)}\cdots Y_{n,\sigma(n)}, 
\] 
where the sum is
taken over the permutations of $\{1, \dots, n\}$ and $l(\sigma)$ is the usual
length function on such permutations. The quantum determinant is a central
element of $\oqmn$, see, for example, \cite[Theorem 4.6.1]{pw}. If $I$ is a
$t$-element subset of $\{1, \dots, m\}$ and $\Gamma$ is a $t$-element subset
of $\{1, \dots, n\}$, then the quantum determinant of the subalgebra of
$\oqmmn$ generated by $Y_{i,\alpha}$, with $i\in I$ and $\alpha \in \Gamma$,
is denoted by $[I\mid \Gamma]$. The elements $[I\mid \Gamma]$ are the {\em
quantum minors} of $\oqmmn$. Note that the quantum minors are
$\ch$-eigenvectors; and so every ideal generated by a family of quantum
minors is $\ch$-invariant.

It follows from \cite[Proposition 4.2]{llr} that there are exactly $m+n-1$
$\ch$-invariant prime ideals in $R$ that have height one. To make  this
result precise, let us introduce some notation.
 
For $1\leq i \leq n+m-1$, let $b_i$ be the quantum minor defined as follows. 
$$b_i:= \left\{
\begin{array}{ll}
\left[1, \dots , i \mid n-i+1 , \dots , n \right] 
& \mbox{ if }1 \leq i \leq m
\\ \left[1, \dots , m \mid n-i+1 , \dots , n+m-i \right] 
& \mbox{ if } m < i
\leq n \\ \left[i-n+1, \dots , m \mid 1 , \dots , m+n-i \right] 
& \mbox{ if }
n < i \leq m+n-1
\end{array}
\right.$$

Note that the $b_i$ with $m\leq i\leq n$ are precisely the $m\times m$ minors
of $\oqmn$ that have consecutive column indices. They are homogeneous of
degree $m$. This fact will be used several times later.

\begin{proposition}[\cite{llr}, Proposition 4.2]
\label{m+n-1}
There are precisely $m+n-1$ height one prime ideals 
that are $\ch$-invariant in $\oqmmn$. 
They are the ideals generated by $b_1,\dots, b_{m+n-1}$. 
\end{proposition}

The quantum minors $b_i$ are normal elements of $R$. Moreover they belong to
the algebra $\overline{R}$ obtained from $R$ by the standard 
deleting derivations
algorithm, see Section \ref{subsectionStrate0}. Indeed, every quantum minor
$b_i$ can be expressed as a product of the canonical generators $T_{\ia}$ of
$\overline{R}$ as follows.

\begin{lemma} \label{normalT}
For $1 \leq i \leq m+n-1$, we have 
$$b_i= \left\{ \begin{array}{ll}
T_{1,n-i+1}T_{2,n-i+2}\dots T_{i,n} & \mbox{ if } 1 \leq i \leq m \\
T_{1,n-i+1}T_{2,n-i+2}\dots T_{m,n+m-i} & \mbox{ if } m < i \leq n \\
T_{i-n+1,1}T_{i-n+2,2}\dots T_{m,m+n-i} & \mbox{ if } n < i \leq n+m-1 \\
\end{array} \right.$$
\end{lemma}

\begin{proof} This lemma is a consequence of \cite[Proposition 5.2.2]{c2}.
\end{proof}

Recall that two elements $a,b$ of $\fract(R)=\fract(\overline{R})$ are said to
$q$-commute if there exists an integer $\bullet$ such that $ab=q^{\bullet}
ba$. Since the $T_{i,\alpha}$ $q$-commute pairwise, it follows from the
previous Lemma \ref{normalT} that the $b_j$ $q$-commute with the
$T_{i,\alpha}$ and that the $b_j$ also $q$-commute pairwise. 
Sometimes, it will not be necessary to know exactly the integers that
appear in the power of $q$. However, at some points, we will need the
following commutation relations that can be easily deduced from Lemma
\ref{normalT} and from the commutation relations between the $T_{i,\alpha}$.

\begin{corollary}\label{commutationbi}
Assume that $m<n$. Then,

\begin{enumerate}

\item for all $m \leq i < j \leq n$, 
$$b_i b_j = q^{\alpha_{ij}} b_j
b_i,$$ where $\alpha_{ij}:= \begin{array}{|c|} \{n-i+1,\dots,m+n-i\} \cap
\{n-j+1,\dots,m+n-j\} \end{array} -m$,  

\item for all $i \in \{1,\dots,m-1\}$,
$$b_n b_i = q^{\alpha_i} b_i b_n \quad 
\mbox{  and } \quad b_n b_{m+n-i} =
b_{m+n-i} b_n,$$ where $\alpha_{i}:= \begin{array}{|c|} \{n-i+1,\dots,n\} \cap
\{m+1,\dots,n\} \end{array}$. 

\end{enumerate}
\end{corollary}

\subsection{The $\ideal{0}$-stratum of $\spec(\oqmmn)$.}
\label{sectionm'n'}

First, it follows from Theorem~\ref{theorem0stratum1}  
that 
the $\ideal{0}$-stratum of
$\spec(\oqmmn)$ is reduced to $\ideal{0}$ when  $v_2(m) \neq v_2(n)$. 
Thus, throughout this section, we
assume that $v_2(m) = v_2(n)$; 
so that $\spec_{\ideal{0}} (\oqmmn)$ is
Zariski-homeomorphic to the prime spectrum of a commutative Laurent polynomial
ring in $m \wedge n$ indeterminates. We set $d:=m \wedge n$ and we denote by
$m'$ and $n'$ the positive integers such that $m=dm'$ and $n=dn'$. 
Observe
that $m'$ and $n'$ are odd, since $v_2(m) = v_2(n)$. 
This observation will be
crucial in what follows.

In Section \ref{subsectionStrate0}, we have shown that $\spec_{\ideal{0}} (R)$
is Zariski-homeomorphic to the prime spectrum of the quantum torus
$P(\Lambda)=\overline{R} \Sigma^{-1}$, where $\Sigma$ denotes the
multiplicative system of $\overline{R}$ generated by the normal elements
$T_{\ia}$ with $(\ia) \in \gc 1,m \dc \times \gc 1,n \dc$. 
Moreover, 
$\spec(P(\Lambda))$ is Zariski-homeomorphic
via extension and contraction to the prime spectrum of the centre
$Z(P(\Lambda))$ of $P(\Lambda)$, 
by \cite[Corollary 1.5]{glettore}. Before describing the $\ideal{0}$-stratum
of $R$, we calculate the centre $Z(P(\Lambda))$ of the quantum torus
$P(\Lambda)$

\subsubsection{The centre of the quantum torus
$P(\Lambda)=\overline{R}\Sigma^{-1}$}

Recall that the quantum minors $b_i$ belong to $\overline{R}$. Moreover, 
the $b_i$ are invertible in
the quantum torus $P(\Lambda)=\overline{R} \Sigma^{-1}$,

For $j \in \{ 1, \dots, d\}$, set 
$$\Delta_j:= \prod_{i=0}^{m'+n'-1}b_{id+j}^{(-1)^i}.$$
(Here we set $b_{m+n}:=1$.)

\begin{theorem} \label{CenterQantumTorus}
$Z(P(\Lambda))= \mathbb{C}[\Delta_1^{\pm 1},\dots , \Delta_d^{\pm 1}]$
\end{theorem}
\begin{proof} First, straightforward computations (see also \cite[Theorem 2.13]{mp}) show that
$\Delta_1,\dots,\Delta_d$ are central in $P(\Lambda)$, so that 
$Z(P(\Lambda))\supseteq \mathbb{C}[\Delta_1^{\pm 1},\dots , \Delta_d^{\pm 1}]$.

If $\underline{s}=(s_{1,1},s_{1,2},\dots,s_{m,n}) \in \mathbb{Z}^{mn}$, we set
$T^{\underline{s}} := T_{1,1}^{s_{1,1}}T_{1,2}^{s_{1,2}} \dots
T_{m,n}^{s_{m,n}} \in P( \Lambda)$. Since the $T_{\ia}$ $q$-commute, we deduce
from Lemma \ref{normalT} that the central elements $\Delta_i$ can be expressed
as follows: $$\Delta_i= q^{\gamma_i} T^{\underline{u^{(i)}}},$$ where
$\gamma_i \in\mathbb{Z}$ and
$\underline{u^{(i)}}=(u_{1,1}^{(i)},u_{1,2}^{(i)},\dots,u_{m,n}^{(i)}) \in
\mathbb{Z}^{mn}$ with $(u_{1,m-d+1}^{(i)}, u_{1,m-d+2}^{(i)},\dots,
u_{1,m}^{(i)})=(0,\dots,0,1,0,\dots,0)$ (the one being in the $i$-th
position).

As in \cite{glettore}, we denote by $\sigma : \mathbb{Z}^{mn} \times
\mathbb{Z}^{mn} \rightarrow \mathbb{C}^*$ the antisymmetric bicharacter
defined by 
$$
\sigma(\underline{s},\underline{t})=\prod_{k,l=1}^{mn}
\Lambda_{k,l}^{s_kt_l}$$
for all $ \underline{s},\underline{t} \in
\mathbb{Z}^{mn}.$ 

Then it follows from \cite[1.3]{glettore} that the centre
$Z(P(\Lambda))$ of $P(\Lambda)$ is a Laurent polynomial ring in the variables
$(T^{\underline{b_1}})^{\pm 1},\dots,(T^{\underline{b_r}})^{\pm 1}$, where
$(\underline{b_1},\dots,\underline{b_r})$ is any basis of $S:=\{\underline{s}
\in \mathbb{Z}^{mn} \mid \sigma(\underline{s},-)\equiv 1\}$.

Now, because of Theorem \ref{theorem0stratum1}, we know $r=\mathrm{rk}(S)=d$.
Moreover, the $\Delta_i$ are central, so that the $\underline{u^{(i)}}$ belong to $S$. To
conclude, observe that, since we have $$(u_{1,m-d+1}^{(i)},
u_{1,m-d+2}^{(i)},\dots, u_{1,m}^{(i)})=(0,\dots,0,1,0,\dots,0),$$
the group  $\mathbb{Z}^{mn}/\sum_{i=1}^{d} \mathbb{Z} \underline{u^{(i)}}$ is torsionfree. 
Hence, the $\underline{u^{(i)}}$ form a basis of $S$, and so the centre $Z(P(\Lambda))$ of $P(\Lambda)$
is a Laurent polynomial ring in the variables $\Delta_1^{\pm
1},\dots,\Delta_d^{\pm 1}$, as desired.
\end{proof}

\subsubsection{Height one primes in $\oqmmn$.}

Let $P$ a height one prime of $R=\oqmmn$. Because of the
$\ch$-stratification, see (\ref{partition}) in Section 1.1, 
there exists an $\ch$-invariant
prime ideal $J$ of $R$ such that $P$ belongs to the $\ch$-stratum associated
to $J$. In particular, we have $J \subseteq P$. Since $P$ has height one, this
implies that the height of $J$ is at most $1$. Thus, two cases arise.
\begin{enumerate}
\item First, suppose that $J$ has height one. In this case, $P= J$, 
since $P$ has height one; 
so that $P$ is an $\ch$-invariant height one prime ideal of $R$. 
Hence it follows from Proposition \ref{m+n-1} that there exists 
$i \in \{1,\dots,m+n-1\}$ such that $P=\ideal{b_i}$.
\item Next, suppose that $J=\ideal{0}$. In this case, $P$ is a height 
one prime ideal that belongs to the $\ideal{0}$-stratum of $\spec(R)$. 
Note that this case can only arise when  $v_2(m)= v_2(n)$, since 
it follows from 
Theorem \ref{theorem0stratum1} that the $\ideal{0}$-stratum of 
$\spec(\oqmmn)$ is reduced to $\ideal{0}$ when 
$v_2(m) \neq v_2(n)$. 
\end{enumerate}
$ $

Let $\spec^1(R)$ denote the set of all height one primes of $R$. The previous
discussion proves the following statement.

\begin{proposition}\label{heightone}$ $
\begin{enumerate}
\item If $v_2(m) \neq v_2(n)$, then $$\spec^1(R)=\left\{ \ideal{b_i} \mid i
\in \left\{1,\dots,m+n-1 \right\} \right\}.$$ \item If $v_2(m) = v_2(n)$, then
$$
\spec^1(R)=\left \{ \ideal{b_i} 
\mid i \in \left\{1,\dots,m+n-1 \right\} \right\}
\cup \spec_{\ideal{0}}^1(R),
$$ 
where $\spec_{\ideal{0}}^1(R)$ denotes the set
of all height one primes of $R$ that belong to the $\ideal{0}$-stratum of
$\spec(R)$.
\end{enumerate}
\end{proposition}

In order to complete the previous result, we now describe, in the case where
$v_2(m) = v_2(n)$, the height one primes of $\oqmmn$ that belongs to
$\spec_{\ideal{0}}(\oqmmn)$.

\begin{proposition} \label{0stratum}
Assume that $v_2(m) = v_2(n)$. Then, for any height one prime ideal $P$ of
$\oqmmn$ that belongs to $\spec_{\ideal{0}}(\oqmmn)$, there exists a
unique (up to scalar) irreducible polynomial $V=\sum_{i_1=0}^{r_1} \dots \sum_{i_d=0}^{r_d}
a_{i_1,\dots,i_d} X_1^{i_1} \dots X_d^{i_d} \in \mathbb{C}[X_1,\dots,X_d]$
(where $r_i = \deg_{X_i}V$) with $V \neq X_i$ 
for all $i \in \{1,\dots,d\}$ such that
$
P = \ideal{u} $
where 
$$
u:=\sum_{i_1=0}^{r_1} \dots
\sum_{i_d=0}^{r_d} a_{i_1,\dots,i_d} \prod_{j=1}^d \left[
\prod_{\substack{i=0\\ i \mbox{ even}}}^{m'+n'-1} b_{id+j} \right]^{i_j}
\left[ \prod_{\substack{i=0\\ i \mbox{ odd}}}^{m'+n'-1}
b_{id+j}\right]^{r_j-i_j}.  
$$ 
Moreover, $u$ is normal in $R$. 
\end{proposition}
\begin{proof} We proceed in three steps.\\$ $

\noindent {\it $\bullet$ Step 1. A generator for the extension of $P$ in a
localisation of $R$.}\\$ $

First, observe that the prime ideals in $\spec_{\ideal{0}} (R)$ do not contain
any $b_i$. Indeed, assume that this is not the case; that is, assume that
there exists $P \in \spec_{\ideal{0}} (R)$ with $b_i \in P$ for a certain $i$.
Then, since $b_i$ is an $\ch$-eigenvector, we have $b_i \in \bigcap_{h\in \ch}
h.P = \ideal{0}$. This is a contradiction; and so 
$\spec_{\ideal{0}} (R) \subseteq \{ P \in \spec(R) \mid b_i \notin P \mbox{
for all } i\}$. On the other hand, if $P$ is a prime ideal of $R$ such that
$b_i \notin P$ for all $i$, then $\bigcap_{h\in \ch} h.P $ is an
$\ch$-invariant prime ideal of $R$ that does not contain any $b_i$. However,
because of \cite[Proposition 2.9]{llr}, every nonzero $\ch$-invariant prime
ideal of $R$ contains a height one prime that is $\ch$-invariant. In other
words, every nonzero $\ch$-invariant prime ideal of $R$ contains a $b_i$. Thus
$\bigcap_{h\in \ch} h.P =\ideal{0}$ and so $P \in \spec_{\ideal{0}} (R)$. To
sum up, we have shown that $\spec_{\ideal{0}} (R) = \{ P \in \spec(R) \mid b_i
\notin P \mbox{ for all } i \}$.

Denote by $T$ the localisation of $R$ with respect of the multiplicatively
closed set $\mathcal{B}$ generated by the normal elements $b_i$. 
Note that 
the torus $\ch$ still acts rationally
by automorphisms on $T$, 
since the $b_i$ are $\ch$-eigenvectors, 
see \cite[Exercise II.3.A]{bg}. 
Moreover, it follows
from the previous study, and from classical results of non-commutative
localisation theory, 
that the map $\varphi : P \rightarrow P\mathcal{B}^{-1}$
is an increasing bijection from $\spec_{\ideal{0}} (R)$ onto $\spec \left(T
\right)$.

Next, 
$T$ is $\ch$-simple; 
that is, the only $\ch$-invariant ideals in $T$ are $\ideal{0}$ and $T$, 
by \cite[Proposition 3.5]{llr}. The 
action of $\ch$ on $T$ is rational, see \cite[Exercise II.3.A]{bg}, and
$T$ is $\ch$-simple; so \cite[Corollary II.3.9]{bg} applies to $T$. Thus
extension and contraction provide mutually inverse bijections between
$\spec(T)$ and $\spec(Z(T))$.

Now, recall that the $b_i$ can be expressed as products of $T_{ \ia}$, see
Lemma \ref {normalT}; so  the $b_i$ belong to the quantum torus
$P(\Lambda)=\overline{R} \Sigma^{-1}$. Moreover, it follows from
\cite[Th\'eor\`eme 3.3.1]{cauchoneff} that there exists a multiplicative
system $S$ of $R$ such that $R \subseteq RS^{-1} = P(\Lambda) = \overline{R}
\Sigma^{-1}$. Hence $T$ is a subalgebra of the quantum torus $P(\Lambda)=R
S^{-1}$, and $R \subseteq T \subseteq P(\Lambda)=\overline{R} \Sigma^{-1}=R
S^{-1}$. Thus $Z(T) \subseteq Z(P(\Lambda)) = \mathbb{C}[\Delta_1^{\pm
1},\dots,\Delta_d^{\pm 1}]$. But the $\Delta_i$ are just products of $b_i^{\pm
1}$, so that they belong to $T$. Hence $Z(T)= \mathbb{C}[\Delta_1^{\pm
1},\dots,\Delta_d^{\pm 1}]$.

Observe that 
$\Delta_j$ can be written as
follows. $$\Delta_j = q^{\bullet} \prod_{\substack{i=0\\ i \mbox{
even}}}^{m'+n'-1} b_{id+j} \prod_{\substack{i=0\\ i \mbox{ odd}}}^{m'+n'-1}
b_{id+j}^{-1},$$ where as usual $\bullet$ denotes an integer,  
since the $b_i$ $q$-commute. For $j \in
\{1,\dots,d\}$, we set
\begin{eqnarray}
\label{defDeltaprimej}
\Delta'_j = \left( \prod_{\substack{i=0\\ i \mbox{ even}}}^{m'+n'-1} b_{id+j}
\right) \left( \prod_{\substack{i=0\\ i \mbox{ odd}}}^{m'+n'-1} b_{id+j}
\right)^{-1}.
\end{eqnarray}
Since $Z(T)= \mathbb{C}[\Delta_1^{\pm 1},\dots,\Delta_d^{\pm 1}]$, the centre
of $T$ is also the (commutative) Laurent polynomial ring in the indeterminates
$\Delta'_1,\dots,\Delta'_d$, that is, 
$$
Z(T)= \mathbb{C}[\Delta_1^{'\pm
1},\dots,\Delta_d^{'\pm 1}].
$$ $ $

Now, let $P$ be a height one prime of $R$ that belongs to $\spec_{\ideal{0}}
(R)$. It follows from the previous study that $\varphi (P) \cap Z(T)$ is
a height one prime of $Z(T)= \mathbb{C}[\Delta_1^{'\pm
1},\dots,\Delta_d^{'\pm 1}]$. 
Hence there exists an irreducible polynomial $V=\sum_{i_1=0}^{r_1} \dots
\sum_{i_d=0}^{r_d} a_{i_1,\dots,i_d} X_1^{i_1} \dots X_d^{i_d} \in
\mathbb{C}[X_1,\dots,X_d]$, with $r_i = \deg_{X_i}V$, and with $V \neq X_i$
for all $i \in \{1,\dots,d\}$, such that 
$$
\varphi (P) \cap Z(T)=\ideal{\;
\sum_{i_1=0}^{r_1} \dots \sum_{i_d=0}^{r_d} a_{i_1,\dots,i_d} \Delta_1^{'i_1}
\dots \Delta_d^{'i_d}}_{Z(T)\;\;}. 
$$ 
Thus, 
$$
\varphi (P)=\ideal{ \;\sum_{i_1=0}^{r_1} \dots \sum_{i_d=0}^{r_d}
a_{i_1,\dots,i_d} \Delta_1^{'i_1} \dots \Delta_d^{'i_d}\;}_T,
$$ 
since extension and contraction provide
mutually inverse bijections between $\spec(T)$ and $\spec(Z(T))$. 
Since the
$b_i$ are invertible in $T$, (\ref{defDeltaprimej}) leads to
$\varphi(P)=\ideal{u}$, where 
$$
u:=\sum_{i_1=0}^{r_1} \dots \sum_{i_d=0}^{r_d} a_{i_1,\dots,i_d}
\prod_{j=1}^d \left[ \prod_{\substack{i=0\\ i \mbox{ even}}}^{m'+n'-1}
b_{id+j} \right]^{i_j} \left[ \prod_{\substack{i=0\\ i \mbox{ odd}}}^{m'+n'-1}
b_{id+j} \right]^{r_j-i_j} .
$$ 
Note that $u$ is a normal element in both $R$
and $T$. \\$ $

\noindent {\it $\bullet$ Step 2. We prove that $u$ is not contained in any 
$\ideal{b_k}$.}\\$ $

Set
$$
u_{i_1,\dots,i_d}:= \prod_{j=1}^d \left[ \prod_{\substack{i=0\\ i \mbox{
even}}}^{m'+n'-1} b_{id+j} \right]^{i_j} \left[ \prod_{\substack{i=0\\ i
\mbox{ odd}}}^{m'+n'-1} b_{id+j} \right]^{r_j-i_j}
$$
for all 
$(i_1,\dots,i_d) \in \gc 0,r_1 \dc \times \dots \times \gc 0,r_d \dc$; 
so that 
$$u=\sum_{i_1=0}^{r_1} \dots \sum_{i_d=0}^{r_d} a_{i_1,\dots,i_d}
u_{i_1,\dots,i_d}.$$
Since the $b_k$ are $\ch$-eigenvectors, it is clear that the
$u_{i_1,\dots,i_d}$ are also $\ch$-eigenvectors.
Hence, for all $(i_1,\dots,i_d) \in \gc 0,r_1 \dc \times \dots \times \gc
0,r_d \dc$, there exists a (unique)
character $f_{i_1,\dots,i_d} : \ch \rightarrow \mathbb{C}^*$ such that
$h.u_{i_1,\dots,i_d} = f_{i_1,\dots,i_d}(h) u_{i_1,\dots,i_d}$ for all $h \in
\ch$. The character $f_{i_1,\dots,i_d}$ is called the $\ch$-eigenvalue of the
$\ch$-eigenvector
$u_{i_1,\dots,i_d}$.

\begin{claim}
\label{claim1}
The $u_{i_1,\dots,i_d}$ are $\ch$-eigenvectors with pairwise distinct
$\ch$-eigenvalues.
\end{claim} 
\noindent {\it Proof of Claim \ref{claim1}.} Let $(i_1,\dots,i_d) \in \gc
0,r_1 \dc \times \dots \times \gc 0,r_d \dc$ and
$h:=(g_1,\dots,g_m;h_1,\dots,h_n) \in \ch$. Then,  
$$
h.u_{i_1,\dots,i_d} = g_1^{\bullet} \dots g_m^{\bullet} h_1^{\bullet} \dots
h_{n-m}^{\bullet}h_{n-m+1}^{\alpha_{m-1}} \dots h_n^{\alpha_0}
u_{i_1,\dots,i_d}, 
$$ 
where $\bullet \in \mathbb{Z}$ and $\alpha_k$ is the
number of quantum minors involving the $(n-k)$-th column in the product that
defines $u_{i_1,\dots,i_d}$.

Now one can observe that, for all $k \in \{0,\dots,d-1\}$, the $(n-k)$-th
column appears in $b_{id+j}$ if and only if either 
(i) $k+1 \leq j \leq d $ and
$0 \leq i \leq m'-1$, or 
(ii) $1 \leq j \leq k$ and $1 \leq i \leq m'$. \\
Hence,
since $m'$ is odd, 
\begin{eqnarray*}
\alpha_k & = & \sum_{j=1}^k \left( \frac{m'-1}{2}i_j + 
\frac{m'+1}{2}(r_j-i_j) \right)
+ \sum_{j=k+1}^d \left( \frac{m'+1}{2}i_j + \frac{m'-1}{2}(r_j-i_j) \right)\\
 & = & \sum_{j=1}^k \left( \frac{m'+1}{2}r_j-i_j \right)
+ \sum_{j=k+1}^d \left( \frac{m'-1}{2}r_j + i_j \right)\\
\end{eqnarray*}

Now, let
$(i_1,\dots,i_d),(i'_1,\dots,i'_d)
\in \gc 0,r_1 \dc \times \dots \times \gc 0,r_d \dc$,
and assume that $u_{i_1,\dots,i_d}$ and $u_{i'_1,\dots,i'_d}$
are associated to the same $\ch$-eigenvalue. Then, since $\mathbb{C}$
is infinite, it follows from the previous study that,
$$
\sum_{j=1}^k \left( \frac{m'+1}{2}r_j-i_j \right)
+ \sum_{j=k+1}^d \left( \frac{m'-1}{2}r_j + i_j \right)
=\sum_{j=1}^k \left( \frac{m'+1}{2}r_j-i'_j \right)
+ \sum_{j=k+1}^d \left( \frac{m'-1}{2}r_j + i'_j \right),
$$
for all $k \in \{0,\dots,d-1\}$.
Hence,
$$
\sum_{j=1}^k \left( i'_j-i_j \right)+ \sum_{j=k+1}^d \left( i_j-i'_j
\right)=0,
$$
for all $k \in \{0,\dots,d-1\}$. 
This forces $i_j=i'_j$ for all $j \in \{ 1,\dots,d \}$; so the claim is
proved. 
$\square$
\\$ $

Next, we 
prove that $u$ does not belong to any $\ideal{b_k}$. Indeed, assume that 
$u \in \ideal{b_k}$ for a certain $k \in \{ 1, \dots, m+n-1\}$. 
Then, since $\ideal{b_k}$ is an $\ch$-invariant prime ideal of $R$, 
it follows from the previous claim and from \cite[II.2.10]{bg} that:
\begin{eqnarray}
\label{critere}
\mbox{If }a_{i_1,\dots,i_d} \neq 0 \mbox{, then } 
u_{i_1,\dots,i_d}= \prod_{j=1}^d \left[ \prod_{\substack{i=0\\ 
i \mbox{ even}}}^{m'+n'-1} b_{id+j} \right]^{i_j} 
\left[ \prod_{\substack{i=0\\ 
i \mbox{ odd}}}^{m'+n'-1} b_{id+j} \right]^{r_j-i_j} \in \ideal{b_k}.
\end{eqnarray}

Set $k:= rd+s$ with $r \in \{0,\dots, m'+n'-1\}$ and $s \in \{1,\dots,d\}$. We
distinguish between two cases. First, assume that $r$ is even. Since the
$\ideal{b_i}$ 
are pairwise distinct height one completely prime ideals of $R$, we
deduce from (\ref{critere}) that $a_{i_1,\dots,i_d} = 0$ if $i_s=0$. This
implies that $V= \sum_{i_1=0}^{r_1} \dots \sum_{i_d=0}^{r_d} a_{i_1,\dots,i_d}
X_1^{i_1} \dots X_d^{i_d} = X_s V'$, with $V'\in \mathbb{C}[X_1,\dots,X_d]$.
This contradicts the facts that $V$ is irreducible and $V \neq X_s$.

Next, assume that $r$ is odd. Since the $\ideal{b_i}$ 
are pairwise distinct height
one completely prime ideals of $R$, we deduce from (\ref{critere}) that
$a_{i_1,\dots,i_d} = 0$ if $i_s=r_s$. Now this contradicts $\deg_{X_s}V=r_s$.

To sum up: $u$ does not belong to any $\ideal{b_k}$.
\\$ $

\noindent {\it $\bullet$ Step 3. We prove that $P$ is generated by $u$.}\\$ $

Recall from step one that $\varphi(P)$ is generated by $u$. 
Hence, it is clear
that $P=\varphi(P) \cap R \supseteq \ideal{u}$. Let now $x \in P$. It
remains to prove that $x \in \ideal{u}$. There exists $(\alpha_1,\dots ,
\alpha_{m+n-1}) \in \mathbb{N}^{m+n-1}$ such that $x b_1^{\alpha_1} \dots
b_{m+n-1}^{\alpha_{m+n-1}} = u r $ with $r \in R$. Choose such a
$(\alpha_1,\dots , \alpha_{m+n-1}) \in \mathbb{N}^{m+n-1}$ minimal (for the
lexicographic order). If $(\alpha_1,\dots , \alpha_{m+n-1}) \neq 0$, there
exists $k$ such that $\alpha_k \neq 0$. Then $ur = x b_1^{\alpha_1} \dots
b_{m+n-1}^{\alpha_{m+n-1}}$ belongs to the completely prime ideal of $R$
generated by $b_k$. 
Hence, $u \in \ideal{b_k}$ or $r \in \ideal{b_k}$. Because of step two,
the first possibility can not happen. 
Hence, $r \in \ideal{b_k}$. Since $b_k$ is
normal, we can write $r= r' b_k$ with $r' \in R$. Thus $x b_1^{\alpha_1} \dots
b_{m+n-1}^{\alpha_{m+n-1}} = u r = u r' b_k $; and so $x b_1^{\alpha_1} \dots
b_k^{\alpha_k -1} \dots b_{m+n-1}^{\alpha_{m+n-1}} = u r' $. This 
contradicts the minimality of $(\alpha_1,\dots , \alpha_{m+n-1})$.
Hence $(\alpha_1,\dots , \alpha_{m+n-1}) = 0$; 
and so $x=ur \in \ideal{u}$, as
desired.
\end{proof}




\section{Automorphisms of quantum matrices}

In this section, we investigate the group of automorphisms of $\oqmmn$. Using
graded arguments together with the results of the previous sections, we show
that, in the non-square case, every $\ch$-invariant height one prime of
$\oqmmn$, except possibly one, is invariant under every automorphism. Next,
by using the preferred basis of $\oqmmn$ introduced in \cite{glduke}, 
we show that the
group of automorphisms of $\oqmmn$, with $2 \leq m<n$, 
is isomorphic to the torus
$(\mathbb{C}^*)^{m+n-1}$.

In the sequel, we will use several times the following 
well-known  result concerning 
normal elements of $R$.

\begin{lemma}
\label{utile}
Let $u$ and $v$ two nonzero normal elements of $R$ such that $\ideal{u}
=\ideal{v}$. Then
there exist $\lambda,\mu \in \mathbb{C}^*$ such that $u=\lambda v $ and $v =
\mu u$.
\end{lemma}

\subsection{$q$-commutation, gradings and automorphisms.}

Let $A=\oplus_{i \in \N} A_i$ be a $\N$-graded $\mathbb{C}$-algebra with
$A_0=\mathbb{C}$. Assume that $A$ is a domain generated as an algebra by
$x_1$,..., $x_n$, and that $A_1=\mathbb{C}x_1\oplus \dots \oplus
\mathbb{C}x_n$. We set $A_{\geq d}:= \oplus_{i \geq d} A_i$. The following
result was inspired by a result in \cite{alevchamarie}.

\begin{proposition} \label{graduation}
Assume that, for all $i \in \{1,\dots,n \}$, there exist $j \neq i$ and
$q_{ij} \neq 1 $ such that $x_i x_j = q_{ij} x_j x_i$. Let $\sigma $ be an
automorphism of $A$ and $x$ be a nonzero homogeneous element of degree $d$ of
$A$. \\Then $\sigma(x)=y_d + y_{> d}$, where $y_d \in A_d \setminus \{0\}$ and
$y_{>d} \in A_{> d}$. 
\end{proposition}

\begin{proof}
First, observe that it is sufficient to prove that 
$\sigma(A_d) \subseteq A_{\geq d}$, for every automorphism
$\sigma $ of $A$. Indeed, assume that this is the case, and let $x$
be a nonzero homogeneous element of degree $d$ of $A$. Then we can write
$\sigma(x)=y_d + y_{> d}$, where $y_d \in A_d $ and $y_{> d} \in A_{> d}$. If
$y_d=0$, then $\sigma(x) = y_{> d} \in A_{> d}$, and thus
$\sigma^{-1}(A_{>d})$ is not contained in $A_{>d}$. This is a contradiction.

Hence it just remains to prove that, for every automorphism $\sigma $ of $A$,
we have $\sigma(A_d) \subseteq A_{\geq d}$. Naturally, it is sufficient to
prove this result when $x=x_i$ is one of the
canonical generators of $A$. So let $\sigma$ be an automorphism $\sigma$ of
$A$ and $i \in \{ 1,\dots,n \}$. We can write 
$$
\sigma(x_i)= \alpha_i + f_i,
$$
where $\alpha_{i} \in \mathbb{C}$ and $f_i \in A_{\geq 1}$. We have to prove
that $\alpha_i =0$.

Now, by hypothesis, there exist $j \in \{1,\dots,n\}$ and $q_{ij}
\neq 1$ such that $x_i x_j = q_{ij} x_j x_i$. Set 
$$
\sigma(x_j)= f_j+g_j, 
$$ 
where $f_j \in A_t$, $f_j \neq 0$ and $g_j \in A_{> t}$. Applying $\sigma$
to the equality $x_i x_j = q_{ij} x_j x_i$, and next identifying the
homogeneous part of degree $t$ yields $ \alpha_i f_j = q_{ij} f_j \alpha_i$.
Thus, $\alpha_i f_j =0$, since $q_{ij} \neq 1$. Now, since $f_j \neq 0$, this
forces $\alpha_i =0$, as desired.
\end{proof}

Note that the commutativity hypothesis of the Proposition
\ref{graduation} is satisfied by the algebra $R=\oqmmn$, provided that $n
\geq 2$. Indeed, the relations that define $R$ are all quadratic, so that
$R=\oplus_{i \in \N} R_i$ is a $\mathbb{N}$-graded algebra, the canonical
generators $Y_{i,\alpha}$ of $R$ having degree one.

Next, for all $(i,\alpha) \in \gc 1,m \dc \times \gc 1,n \dc$ and $\beta \neq
\alpha$, we have $Y_{i,\alpha} Y_{i,\beta} = q^{\bullet} Y_{i,\beta}
Y_{i,\alpha}$ with $\bullet=\pm 1$. Thus, the commutativity hypothesis of
Proposition \ref{graduation} is satisfied; 
and so one can state:

\begin{corollary}
 \label{graduation2} Let $\sigma$ be an automorphism of $R=\oqmmn$ and $x$ an
 homogeneous element of degree $d$ of $R$. \\Then $\sigma(x)=y_d + y_{> d}$,
 where $y_d \in R_d \setminus \{0\}$ and $y_{>d} \in R_{> d}$. 
\end{corollary}

 Note, for later use, that a
$t \times t$ quantum minor of $R$ is a homogeneous element of degree $t$
with respect to this grading of $R$. In the sequel, $R$ will always be
endowed with this grading.

\subsection{The 
action of $\aut(\oqmmn)$ on the set of height one primes: the
non-square case.}

Throughout this section, we assume that $m < n$. Our aim in this section is to
show that every $\ch$-invariant height one prime of $\oqmmn$, except possibly
one, is invariant under every automorphism. In order to do this, 
we will distinguish
between two cases.

\subsubsection{The case where $n \neq 3m$.}

Throughout this section, we assume that $m < n$ and $n \neq 3m$. In
this case, we first show that the set of all $\ch$-invariant height one primes of $R$ is invariant under every automorphism of $R$, that is, we have:

\begin{lemma}\label{imagebi1}
Assume that $n \neq 3m$. Let $\sigma$ be an automorphism of $R$ and $i \in
\{1,\dots,m+n-1\}$. Then there exists $j \in \{1,\dots,m+n-1\}$ such that
$\sigma( \ideal{b_i}) = \ideal{b_j} $.
\end{lemma}

\begin{proof} First, assume that 
$v_2(m) \neq v_2(n)$. It follows from Proposition
\ref{heightone} that $\{ \ideal{b_i} \mid i \in \gc 1,m+n-1 \dc \}$ is exactly the
set of all height one primes of $R$ and so the result is obvious in this case.

Next, assume 
that $v_2(m) = v_2(n)$. We use the notation of Section
\ref{sectionm'n'}. In particular, $d$ denotes the greatest common divisor of
$m$ and $n$, and we set $m=dm'$ and $n=dn'$. Observe that $m'$ and $n'$ are
odd. Note further that, since $n \neq 3m$, we have $(m',n') \neq (1,3)$.
Moreover, since $m < n$, we have $m' < n'$. Since $m'$ and $n'$ are both 
odd, this forces  $n' \geq 5$.

Since $\ideal{b_i}$ 
is a height one prime ideal of $R$, its image under 
$\sigma$ is
also a height one prime ideal of $R$. We distinguish between two cases.

First, if $\ideal{\sigma(b_i)}$ is $\ch$-invariant, then it follows from Proposition
\ref{m+n-1} that $\ideal{\sigma(b_i)} = \ideal{b_j}$ 
for some $j$ and so the proof
is complete in this case.

Next, assume that $\ideal{\sigma(b_i)}$ is not $\ch$-invariant. In this case,
$\ideal{\sigma(b_i)}$ is a height one (completely) prime ideal of $R$ which
is not $\ch$-invariant. Thus, $\ideal{\sigma(b_i)}$ belongs to the
$\ideal{0}$-stratum of the prime spectrum of $R$, by Proposition
\ref{heightone}; and so we deduce from Proposition \ref{0stratum} that there
exist $(r_1,\dots,r_d) \in \mathbb{N}^d \setminus \{0\}$ and scalars
$a_{i_1,\dots,i_d} \in \mathbb{C}$ such that 
$
\ideal{\sigma(b_i)}=\ideal{u}$, where 
$$
u:=\sum_{i_1=0}^{r_1} \dots \sum_{i_d=0}^{r_d} a_{i_1,\dots,i_d} \prod_{j=1}^d
\left[ \prod_{\substack{i=0\\ i \mbox{ even}}}^{m'+n'-1} b_{id+j}
\right]^{i_j} \left[ \prod_{\substack{i=0\\ i \mbox{ odd}}}^{m'+n'-1} b_{id+j}
\right] ^{r_j-i_j}.
$$ 
It follows
from Proposition \ref{0stratum} that $u$ is a normal element of $R$. On the
other hand, $\sigma(b_i)$ is normal, since $b_i$ is normal. Thus, 
we deduce from
Lemma \ref{utile} that $u \in \mathbb{C}^* \sigma(b_i)$.

Now $b_i$ is a quantum minor of $R=\oqmmn$, and so $b_i$ is a homogeneous
element of degree less than or equal to $m$. So Corollary \ref{graduation2}
implies that
\begin{eqnarray}
\label{contradiction1}
u  \not\in  R_{>m}.
\end{eqnarray}

On the other hand, since $(r_1,\dots,r_d) \neq 0$, there exists $k$ such that
$r_k \geq 1$. Now, recalling that $m'$ and $n'$ are odd
and that $n' \geq 5$, we have: 
$$ 
\mathrm{deg} \left( \prod_{j=1}^d \left[
\prod_{\substack{i=0\\ i \mbox{ even}}}^{m'+n'-1} b_{id+j} \right]^{i_j}
\left[ \prod_{\substack{i=0\\ i \mbox{ odd}}}^{m'+n'-1}
b_{id+j}\right]^{r_j-i_j} \right) \geq \left\{ \begin{array}{ll}
\mathrm{deg} \left( b_k b_{(m'+1)d+k} \right) & \mbox{ if } i_k \geq 1 \\
\mathrm{deg} \left( b_{m'd+k} b_{(m'+n'-3)d+k} \right) & \mbox{ if } i_k =
0 \\
\end{array}
\right.
$$
Since $m < n$, we have $m' < n'$. Further, $m'$ and $n'$ are odd, and so $m'+1
< n'$.
Hence $m=m'd \leq m'd+k \leq (m'+1)d+k \leq n=n'd$. This implies that
$b_{(m'+1)d+k}$ and $b_{m'd+k}$ are $m \times m$ quantum minors. Thus, in both
cases, we get
$$
\mathrm{deg} \left( \prod_{j=1}^d \left[ \prod_{\substack{i=0\\ i \mbox{
even}}}^{m'+n'-1} b_{id+j} \right]^{i_j} \left[\prod_{\substack{i=0\\ i \mbox{
odd}}}^{m'+n'-1} b_{id+j}\right]^{r_j-i_j} \right)
\geq m+1.
$$
Thus $u$ is a linear combination of terms of degree greater than $m$. This
implies that $u \in R_{>m}$, contradicting
(\ref{contradiction1}); and so the proof is complete.
\end{proof}

In fact, more is true: each $\ch$-invariant height one prime is left invariant by any automorphism, as
the following result shows.

\begin{proposition}\label{imagebi2}
Assume that $n \neq 3m$ and let $\sigma$ be an automorphism of $R$. Then, for
each $i \in \{1,\dots,m+n-1\}$, there exists $\lambda_i \in \mathbb{C}^*$ such
that $\sigma(b_i) = \lambda_i b_i $. 
\end{proposition}

\begin{proof}
Let $i \in \{1,\dots,m+n-1\}$. It follows from Lemma \ref{imagebi1} that there
exists $j \in \{1, \dots, m+n-1\}$ such that $\sigma(\ideal{b_i} ) 
= (\ideal{b_j}) $. Since
$\sigma(b_i)$ and $b_j$ are normal, it follows from Lemma \ref{utile} that
there exists $\lambda_j \in \mathbb{C}^*$ such that $\sigma(b_i)= \lambda_i
b_j$.

Thus, there exist scalars $\lambda_1,\dots,\lambda_{m+n-1} \in \mathbb{C}^*$
and a permutation $s \in S_{m+n-1}$ such that $\sigma(b_i) = \lambda_i
b_{s(i)} $ for all $i \in \{1,\dots,m+n-1\}$.

We will now prove that $s$ is the identity, that is, $s(i) = i $ for all $i
\in \{1,\dots,m+n-1\}$.

First, let $i \in \{m,\dots,n\}$, so that $b_i$ is a $m \times m$ quantum
minor. Then $b_i \in R_{m}$ and it follows from Corollary \ref{graduation2}
that $\sigma(b_i) \in R_{\geq m}$. Hence $b_{s(i)} \in R_{\geq m}$. This
implies that $b_{s(i)} $ is also a $m \times m$ quantum minor, so that $s(i)
\in \{m,\dots,n\}$. Thus $s$ induces a permutation of $\{m,\dots,n\}$.

We now prove with the help of a decreasing induction that $s(i) = i$ for all
$i \in \{m,\dots,n\}$.

It follows from Corollary \ref{commutationbi} that $b_n b_j = \qdot b_j b_n$
with $\bullet \geq 0$ for all $j \in \{m,\dots,n\}$. Then, applying 
$\sigma$ leads to $b_{s(n)} b_{s(j)} = \qdot b_{s(j)} b_{s(n)}$ with $\bullet
\geq 0$ for all $j \in \{m,\dots,n\}$. Since $s$ is a permutation of
$\{m,\dots,n\}$, this implies that $b_{s(n)} b_{j} = \qdot b_{j} b_{s(n)}$
with $\bullet \geq 0$ for all $j \in \{m,\dots,n\}$. Now, if $s(n) \neq n$,
then we have $s(n)+1 \in \{m,\dots,n\}$ and $b_{s(n)} b_{s(n)+1} = q^{-1}
b_{s(n)+1} b_{s(n)}$. This is a contradiction and so $s(n)=n$.

We now assume that $m\leq i < n$. It follows from the induction hypothesis
that $s$ induces a permutation of $\{ m, \dots,i\}$; so that $s(i) \leq i$. By
using a similar argument to that in the previous paragraph, we obtain 
$s(i)=i$.

Hence, $s(i) = i$ for all  $i \in \{m,\dots,n\}$.

Now let $i \in \{ 1, \dots, m-1 \}$. Then $b_i \in R_i$; 
and so it follows from
Proposition \ref{graduation} that $\sigma(b_i) \in R_{\geq i} \setminus R_{ >
i}$. Hence $b_{s(i)} \in R_{\geq i} \setminus R_{ > i}$ so that $b_{s(i)} $ is
also a $i \times i$ quantum minor. This implies that either $s(i)=i$ or $s(i)
= m+n-i$. Note that similar arguments show that either $s(m+n-i)=i$ or
$s(m+n-i) = m+n-i$; so  $s$ induces a permutation of $\{ i, m+n-i\}$.
Observe that it follows from Corollary \ref{commutationbi} that $b_n b_{m+n-i}
= b_{m+n-i} b_n$ and $b_n b_i = \qdot b_i b_n$ with $\bullet > 0$. Since we
have already proved that $s(n)=n$, applying  $\sigma$ leads to $b_n
b_{s(m+n-i)} = b_{s(m+n-i)} b_n$ and $b_n b_{s(i)} = \qdot b_{s(i)} b_n$ with
$\bullet > 0$. This forces $s(i)=i$ and $s(m+n-i)=m+n-i$.
\end{proof}

\subsubsection{The case where $n= 3m$.}

Throughout this section, we assume that $n = 3m$. In this case, we are not
able to prove directly that the set of all $\ch$-invariant height one primes
is invariant under every automorphism of $R=\oqmmn$. However, by 
using arguments similar to  
those developed in the proof of Lemma \ref{imagebi1}, one can
establish the following weaker result.

\begin{lemma}\label{imagebi1bis}
Assume that $n = 3m$. Let $\sigma$ be an automorphism of $R$ and $i \in \{1,\dots,m+n-1\}$. 
\\Then, either 
\begin{enumerate}
\item there exists $j \in \{1,\dots,m+n-1\}$ such that 
$\sigma(\ideal{b_i} ) = \ideal{b_j}
$, 
\\  \noindent or,  
\item there exist $\lambda \in \mathbb{C}^*$ and $\mu \in \mathbb{C}$
such that $\sigma(b_i) = \lambda b_{2m} + \mu b_{m}b_{3m}$.
\end{enumerate}
\end{lemma}

\begin{proof} 
Note that $v_2(m) = v_2(n)$, since $n=3m$. Also, observe  that the
greatest common divisor $d$ of $m$ and $n$ is equal to $m$ and that, if we set
$m=dm'$ and $n=dn'$, then $m'=1$ and $n'=3$.

Since $\ideal{b_i}$ is a height one prime ideal of $R$, its image under 
$\sigma$ is also a height one prime ideal of $R$. We distinguish between two
cases.

If $\ideal{\sigma(b_i)}$ is $\ch$-invariant, then it follows from Proposition
\ref{m+n-1} that $\sigma(\ideal{b_i}) = \ideal{b_j} $ for some 
$j$ and so the proof
is complete in this case.

Assume now that $\ideal{\sigma(b_i)}$ is not $\ch$-invariant. In this case,
$\ideal{\sigma(b_i)}$ is a height one (completely) prime ideal of $R$ which
is not $\ch$-invariant. Thus it follows from Proposition \ref{heightone} that
$\ideal{\sigma(b_i)}$ belongs to the $\ideal{0}$-stratum of the prime spectrum
of $R$; and so, recalling that $d=m$, 
we deduce from Proposition \ref{0stratum} that there exist
$(r_1,\dots,r_m) \in \mathbb{N}^m \setminus \{0\}$ and scalars
$a_{i_1,\dots,i_m} \in \mathbb{C}$ such that
$
\ideal{\sigma(b_i)}= \ideal{u}$ where 
$$
u:=\sum_{i_1=0}^{r_1} \dots \sum_{i_m=0}^{r_m}
a_{i_1,\dots,i_m} \prod_{j=1}^m \left[ \prod_{\substack{i=0\\ i \mbox{
even}}}^{3} b_{im+j} \right]^{i_j} \left[ \prod_{\substack{i=0\\ i \mbox{
odd}}}^{3} b_{im+j} \right] ^{r_j-i_j}.
$$ 
It follows 
that $u$ is a normal element of $R$, by  Proposition \ref{0stratum}.
On the other hand, 
$\sigma(b_i)$ is normal in $R$, since $b_i$ is normal in $R$. .
Hence, 
$u \in \mathbb{C}^*
\sigma(b_i)$, by  Lemma \ref{utile}.

Now $b_i$ is a quantum minor of $R=\oqmmn$; and so $b_i$ is a homogeneous
element of degree less than or equal to $m$. Thus, 
Corollary \ref{graduation2}
implies that
\begin{eqnarray}
\label{contradiction2}
u   \not\in  R_{>m}.
\end{eqnarray}

On the other hand, since $(r_1,\dots,r_m) \neq 0$, there exists $k$ 
such that $r_k \geq 1$. 
We consider three separate cases.

$\bullet$ First, suppose that $k <m$, then 
$$
\mathrm{deg} \left(
\prod_{j=1}^m \left[ \prod_{\substack{i=0\\ i \mbox{ even}}}^{3} b_{im+j}
\right]^{i_j} \left[ \prod_{\substack{i=0\\ i \mbox{ odd}}}^{3}
b_{im+j}\right]^{r_j-i_j} \right) \geq \left\{ \begin{array}{ll}
\mathrm{deg} \left( b_k b_{2m+k} \right) & \mbox{ if } i_k \geq 1 \\
\mathrm{deg} \left( b_{m+k} b_{3m+k} \right) & \mbox{ if } i_k = 0
\end{array}
\right.
$$

Now, $m+k$ and $2m +k$ are both between $m$ and $n$, since $k<m$ and $n=3m$.
Thus, $b_{m+k}$ and $b_{2m+k}$ are $m \times m$
quantum minors. Also, $b_k$ and $b_{3m+k}$ are homogeneous of degree
greater than or equal to $1$. Thus, in both cases, we get 
$$
\mathrm{deg}
\left( \prod_{j=1}^m \left[ \prod_{\substack{i=0\\ i \mbox{ even}}}^{3}
b_{im+j} \right]^{i_j} \left[\prod_{\substack{i=0\\ i \mbox{ odd}}}^{3}
b_{im+j}\right]^{r_j-i_j} \right) \geq m+1.
$$ 
Hence $u$ is a linear
combination of terms of degree greater than $m$. This implies that $u \in
R_{>m}$. This contradicts (\ref{contradiction2}). \\$ $

$\bullet$ Next, assume that $k=m$ and $r_m \geq 2$. In this 
case, one can prove, by 
using similar arguments, that $u$ is a linear combination
of terms of degree greater than $m$; so this case  cannot happen. \\$ $

$\bullet$ Finally, 
assume that $k=m$ and $(r_1,\dots,r_m)=(0,\dots,0,1)$. In this case, 
there exist $\lambda', \mu' \in \mathbb{C}$ such that  
$$u= \lambda'  b_{2m} + \mu' b_{m}b_{3m}. $$
If $\lambda' =0$, then once again $u$ is a linear combination of terms of 
degree greater than $m$, contradicting (\ref{contradiction2}). 
Hence $\lambda' \neq 0$. Since we have already proved that 
$u \in \mathbb{C}^* \sigma (b_i)$, we see that there exist 
$\lambda \in \mathbb{C}^*$ and $\mu  \in \mathbb{C}$ 
such that $\sigma(b_i) = \lambda  b_{2m} + \mu b_{m}b_{3m}$, as desired.
\end{proof}

The following commutation relations can be easily deduced from Corollary
\ref{commutationbi}.
\begin{lemma}\label{b+bb}
\label{relationBI}
For all $m \leq i \leq 3m$: 
$$b_i (\lambda b_{2m} + \mu b_m b_{3m} )
= \left\{\begin{array}{ll} q^{\alpha_{i,2m}} (\lambda b_{2m} + \mu b_m b_{3m}
) b_i & \mbox{ if } i < 2m \\ q^{-\alpha_{i,2m}} (\lambda b_{2m} + \mu b_m
b_{3m} ) b_i & \mbox{ if } i > 2m, \\
\end{array}
\right.
$$
where 
$\alpha_{i,2m}:= \begin{array}{|c|} \{3m-i+1,\dots,4m-i\} \cap
\{m+1,\dots,2m\} \end{array} -m$.
\end{lemma}

We can now obtain the analogous (but slightly weaker) result to
Proposition \ref{imagebi2} in the $n=3m$ case.

\begin{proposition}\label{imagebi2bis}
Assume that $n=3m$ and let $\sigma$ be an automorphism of $R$. Then, for all
$i \in \{1,\dots,4m-1\}$, $i \neq 2m$, there exists $\lambda_i \in
\mathbb{C}^*$ such that $\sigma(b_i) = \lambda_i b_i $.
\end{proposition}

\begin{proof}
Let $i \in \{1,\dots,4m-1\}$ and assume that there exists $j \in \{1, \dots,
4m-1\}$ such that $\sigma(\ideal{b_i}) = \ideal{b_j}$. 
Then, 
it follows from Lemma \ref{utile} that there exists $\lambda_j \in
\mathbb{C}^*$ such that $\sigma(b_i)= \lambda_i b_j$, since $\sigma(b_i)$ and
$b_j$ are normal.

Thus,
we deduce from Lemma \ref{imagebi1bis} that,
for all $i \in \{ 1,\dots,4m-1\}$, either there exist $j \in \{1, \dots,
4m-1\}$ and $\lambda_i\in \mathbb{C}^*$ such that $\sigma(b_i)=\lambda_i b_j$,
or there exist $\lambda_i,\mu_i \in \mathbb{C}$ with $\lambda_i \neq 0$ such
that $\sigma(b_i) = \lambda_i b_{2m} + \mu_i b_m b_{3m} $.

We distinguish between two cases.\\$ $

\noindent $\bullet$ If, for each $i \in \{ 1,\dots,4m-1\}$, there exist $j
\in \{1, \dots, 4m-1\}$ and $\lambda_i\in \mathbb{C}^*$ such that
$\sigma(b_i)=\lambda_i b_j$, then, by 
using similar arguments to those in  the proof
of Proposition \ref{imagebi2}, we show that $\sigma(b_i)=\lambda_i b_i$ for
all $i \in \{1,\dots,4m-1\}$. \\$ $

\noindent $\bullet$ Now, assume that there exist $k \in \{ 1,\dots,4m-1\}$,
and 
$\lambda_k,\mu_k \in \mathbb{C}$ with $\lambda_k \neq 0$, such that
$\sigma(b_k) = \lambda_k b_{2m} + \mu_k b_m b_{3m} $.\\$ $

\noindent $\bullet\bullet$ First, we show that $k=2m$.

Observe that, for all $i \in \{1,\dots,m-1\} \cup \{3m+1,\dots,4m-1\}$, the
quantum minor $b_i$ is a homogeneous element of degree less than $m$. 
Hence,
$\sigma(b_i) \notin R_{\geq m}$, by  Corollary \ref{graduation2}. 
Since $\sigma(b_k) = \lambda_k b_{2m} + \mu_k b_m b_{3m} \in 
R_{\geq m}$, this shows that $i \neq k$. Thus, $k \in \{m,\dots,3m\}$.

Next, let $i \in \{m,\dots,3m\}$, so that $b_i$ is a $m \times m$ quantum
minor. Then $b_i \in R_{ m}$ and it follows from Corollary \ref{graduation2}
that $\sigma(b_i) \in R_{\geq m}$. Hence, either 
\begin{eqnarray}
\label{cas1}
\mbox{there exist $j \in \{m, \dots, 3m\}$ and $\lambda_i\in
\mathbb{C}^*$ such that $\sigma(b_i)=\lambda_i b_j$,} & &
\end{eqnarray}
or 
\begin{eqnarray}
\label{cas2}
\mbox{there exist $\lambda_i,\mu_i \in \mathbb{C}$ with $\lambda_i \neq 0$
such that $\sigma(b_i) = \lambda_i b_{2m} + \mu_i b_m b_{3m} $.} & &
\end{eqnarray}

We now prove that
(\ref{cas2}) can not happen, for all $i \in \{ 2m+1,\dots,3m\}$. 
Indeed, assume that there exists $i \in \{ 2m+1,\dots,3m\}$ such that 
$\sigma(b_i) = \lambda_i b_{2m} + \mu_i b_m b_{3m} $. 
Note that 
$b_i b_j =\qdot b_j b_i $ with $\bullet > 0$, for all $m \leq j <i$, 
by Corollary \ref{commutationbi}.

Hence, by applying $\sigma$, we obtain:
\begin{eqnarray}
\label{equationProposition39}
(\lambda_i b_{2m} + \mu_i b_m b_{3m}) \sigma (b_j) & = & \qdot \sigma( b_j )
(\lambda_i b_{2m} + \mu_i b_m b_{3m})
\end{eqnarray}
with $\bullet > 0$ for all $m \leq j <i$. 
This implies that $\sigma (b_j)$ can
not be equal to $\lambda b_{2m} + \mu b_m b_{3m}$ with $\lambda,\mu \in
\mathbb{C}$, since 
$\lambda b_{2m}
+ \mu b_m b_{3m}$ commutes with $\lambda_i b_{2m} + \mu_i b_m b_{3m}$
by Lemma~\ref{b+bb}. 
Hence
we deduce from (\ref{cas1}) and (\ref{cas2}) that, for each $j$ such that 
$m \leq j <i$, 
there
exist $l \in \{ m, \dots,3m-1\}$ and $\lambda_j \in \mathbb{C}^*$ such that
$\sigma(b_j) = \lambda_j b_l$. Moreover, we deduce from 
(\ref{equationProposition39}) and Lemma~\ref{b+bb} that 
we must have $l \in \{ m, \dots,2m-1\}$. Thus, since $i > 2m$, 
there exist $j \neq j'$
with $m \leq j,j' <i$ with 
$\sigma (b_j) \in \mathbb{C}^* \sigma(b_{j'})$. This is  impossible
since $\sigma$ is an automorphism and
$\ideal{b_j} \neq \ideal{b_{j'}}$, by Proposition~\ref{m+n-1}.

Hence, for all $i \in \{ 2m+1,\dots,3m\}$, there exist $j \in \{m, \dots,
3m\}$ and $\lambda_i\in \mathbb{C}^*$ such that $\sigma(b_i)=\lambda_i b_j$.
In other words, $k \neq i$ for all $i \in \{ 2m+1,\dots,3m\}$.

A similar argument shows that, for all $i \in \{ m,\dots,2m-1\}$, there exist
$j \in \{m, \dots, 3m\}$ and $\lambda_i\in \mathbb{C}^*$ such that
$\sigma(b_i)=\lambda_i b_j$, and so $k \neq i$ for all $i \in \{
m,\dots,2m-1\}$.

In conclusion, the only possibility is $k=2m$. Hence, we have already proved
that:
\begin{enumerate}

\item There exist $\lambda, \mu \in \mathbb{C}$, with 
$\lambda \neq 0$, such that
$\sigma(b_{2m})= \lambda b_{2m} + \mu b_m b_{3m}$. 

\item For all $i \in
\{m,\dots,3m\}$ with $i \neq 2m$, there exist $j \in \{m, \dots, 3m\}$ and
$\lambda_i\in \mathbb{C}^*$ such that $\sigma(b_i)=\lambda_i b_j$. 

\item For
all $i \in \{1,\dots,m-1\} \cup \{3m+1,\dots,4m-1\}$, there exist $j \in \{1,
\dots, 4m-1\}$ and $\lambda_i\in \mathbb{C}^*$ such that
$\sigma(b_i)=\lambda_i b_j$.
\end{enumerate}
$ $

\noindent $\bullet\bullet$ 
We prove by induction that 
$\sigma(b_i)=\lambda_i b_i$, for all $i \in \{
2m+1,\dots,3m\}$.

First, we know that there exist $\lambda_{2m+1} \in \mathbb{C}^*$ and $j \in
\{m, \dots, 3m\}$ such that $\sigma(b_{2m+1})=\lambda_{2m+1} b_j$. It follows
from Corollary \ref{commutationbi} that $b_{2m}b_{2m+1} = q b_{2m+1}b_{2m}$.
Hence, applying $\sigma$ yields: $$(\lambda b_{2m} + \mu b_m b_{3m}) b_j=
q b_j (\lambda b_{2m} + \mu b_m b_{3m}).$$ In view of Lemma
\ref{relationBI}, this forces $j=2m+1$, as desired.

Next, let $i \in \{ 2m+2,\dots,3m\}$. It follows from the previous study that
there exist $j \in \{m, \dots, 3m\}$ and $\lambda_i\in \mathbb{C}^*$ such that
$\sigma(b_i)=\lambda_i b_j$. Moreover we deduce from the induction hypothesis
that $\sigma(b_{i-1}) = \lambda_{i-1} b_{i-1}$. Now, Corollary
\ref{commutationbi} shows that $b_i b_{i-1} = q b_{i-1} b_i$. Applying 
$\sigma$ yields $$b_j b_{i-1} = q b_{i-1} b_j.$$ In view of 
Corollary \ref{commutationbi}, 
this implies that $j=i$, as desired.

Hence, 
$\sigma(b_i)=\lambda_i b_i$, for all $i \in \{ 2m+1,\dots,3m\}$. 
A similar argument shows that 
$\sigma(b_i)=\lambda_i b_i$,  for all
$i \in \{ m,\dots,2m-1\}$. \\$ $

\noindent $\bullet \bullet$ Now, let $i \in \{ 1, \dots, m-1 \}$. Then $b_i
\in R_i$; and so it follows from Proposition \ref{graduation} that $\sigma(b_i)
\in R_{\geq i} \setminus R_{ > i}$. Now, 
there
exist $j \in \{1, \dots, 4m-1\}$ and $\lambda_i\in \mathbb{C}^*$ such that
$\sigma(b_i)=\lambda_i b_j$
because of the previous study. 
Hence we have $b_{j} \in R_{\geq i} \setminus R_{
> i}$, so that $b_{j} $ is also a $i \times i$ quantum minor. This implies
that either $j=i$ or $j = 4m-i$. Now, it follows from Corollary
\ref{commutationbi} that $b_{3m} b_i = \qdot b_i b_{3m}$ with $\bullet > 0$.
Since we have already proved that $\sigma(b_{3m}) = \lambda_{3m} b_{3m}$,
composing by $\sigma$ leads to $b_{3m} b_{j} = \qdot b_{j} b_{3m}$ with
$\bullet > 0$. On the other hand, it follows Corollary \ref{commutationbi}
that $b_{3m} b_{4m-i} = b_{4m-i} b_{3m}$. So $j$ can not be equal to $4m-i$.
Hence $j=i$ and so $\sigma(b_i)=\lambda_i b_i$, as desired.

A similar argument shows that
there exist $\lambda_i \in \mathbb{C}^*$ such 
that $\sigma(b_i)=\lambda_i b_i$, for all $i \in \{ 3m+1, \dots, 4m-1 \}$.
\end{proof}

\subsection{The automorphism group of non-square quantum matrices.}

\begin{theorem}\label{automxp}
Assume that $m < n$ and $(m,n) \neq (1,3)$. Let $\sigma$ be an automorphism of
$\oqmmn$. Then
there exist $\mu_{i,\alpha} \in \mathbb{C}^*$ such that
$\sigma(Y_{i,\alpha})= \mu_{i,\alpha}Y_{i,\alpha}$, 
for all $(i,\alpha) \in \gc 1, m \dc \times \gc 1,n \dc$.
\end{theorem}

\begin{proof}
We proceed by induction on $n$. The case $n=2$ easily follows from Proposition
\ref{imagebi2}. So we assume that $n \geq 3$. If $m=1$, then $n >3$, and so
once again the result easily follows from Proposition \ref{imagebi2}. So we
assume that $m \geq 2$. We need to see that $\sigma$ acts on each generator
$Y_{i,\alpha}$ by multiplication by a scalar. We do this by using preferred
basis arguments. We use the language and notation of \cite{glduke}. \\$ $

Note that, because of \cite[Proposition 5.3]{glduke}, the quantum minor
$b_n=[1\dots m| 1 \dots m]$ commutes with $Y_{i,\alpha}$, for 
$1 \leq i,\alpha
\leq m$, and that $b_n [R \mid C] = \qdot [R \mid C]b_n$ with $\bullet \geq 0$ for all other quantum minors $[R\mid C]$ of
$\oqmmn$. Thus, $b_n$ $\qdot$-commutes with each monomial in the
preferred basis where $\bullet$ is $\geq 0$, and \cite[Proposition
5.3]{glduke} shows that $\bullet$ is equal to zero if and only if the only
quantum minors that occur in the monomial are those of
$(Y_{i,\alpha})_{i,\alpha \in \{ 1 , \dots, m \}}$. 
Observe further that, if
$b_n=[1\dots m|1 \dots m]$ commutes with an element $y \in \oqmmn$, 
then it must commute
with each of the monomials in the expression for $y$ in the preferred basis.
\\$ $

Now, let $\sigma$ be an automorphism of $\oqmmn$. We will show first that
$\sigma(Y_{i,\alpha}) = \mu_{i,\alpha} Y_{i,\alpha}$ for all $1 \leq i,\alpha
\leq m$. \\

Set $y_{i,\alpha}:= \sigma(Y_{i,\alpha})$. Since $Y_{i,\alpha}$ commutes with
$b_n$ and since $\sigma$ acts on $b_n$
by multiplication by a scalar (see Proposition \ref{imagebi2} if $n \neq 3m$
or Proposition \ref{imagebi2bis} if $n=3m$), then $y_{i,\alpha}$ must commute
with $b_n$. Thus, any monomial $z$ in the expression for
$y_{i,\alpha}$ in the preferred basis must also commute with $b_n$. 
This means that the only quantum minors that can occur in $z$ are
those of $(Y_{i,\alpha})_{i,\alpha \in \{ 1 , \dots, m \}}$. In particular,
$\sigma(Y_{i,\alpha}) $ belongs to the subalgebra $R_{m,m}$ of $R=\oqmmn$
generated by $Y_{i,\alpha}$, $1 \leq i,\alpha \leq m$ (which is a copy of
$O_q(M_{m,m})$). Hence $\sigma$ induces an automorphism of $R_{m,m}$. \\$ $

Moreover,  
$\sigma$ acts on the quantum
minors $b_{n+i}=[i+1,\dots,m | 1, \dots, m-i ]$ by multiplication by 
scalars, by  Proposition \ref{imagebi2} (if $n \neq 3m$) or
Proposition \ref{imagebi2bis} (if $n = 3m$). 
In
particular, there exists $\mu_{m,1} \in \mathbb{C}^*$ such that $\sigma
(Y_{m,1})= \mu_{m,1} Y_{m,1}$. 
Let $ j \in \{ 1, \dots, m-1\}$. Then 
$Y_{j,1} Y_{m,1}=q Y_{m,1} Y_{j,1}$; so that $\sigma(Y_{j,1}) Y_{m,1}=q
Y_{m,1} \sigma(Y_{j,1})$. 
Write $\sigma(Y_{j,1})$ in the PBW basis of
$R_{m,m}$: 
$$
\sigma(Y_{j,1}) = \sum_{\underline{\gamma} \in \Gamma}
c_{\underline{\gamma}} Y_{1,1}^{\gamma_{1,1}} Y_{1,2}^{\gamma_{1,2}} \dots
Y_{m,m}^{\gamma_{m,m}},
$$ 
where $\Gamma$ is a finite subset of
$\mathbb{N}^{m^2}$. 
Hence $\sigma(Y_{j,1}) Y_{m,1}=q Y_{m,1} \sigma(Y_{j,1})$
implies that, for all $\underline{\gamma} \in \Gamma$ such that
$c_{\underline{\gamma}} \neq 0$, we have $\sum_{i=1}^{m-1} \gamma_{i,1} -
\sum_{i=2}^{m} \gamma_{m,i}= 1$. Thus, for all $\underline{\gamma} \in \Gamma$
such that $c_{\underline{\gamma}} \neq 0$, there exists $ i \in \{ 1 , \dots ,
m-1 \}$ with $\gamma_{i,1} \neq 0$. 
Denote by $J$,  the prime
ideal of $R$ generated by the $Y_{i,1}$, with $i \in \{ 1,
\dots,m \}$, and, similarly, denote by $J'$ the prime
ideal of $R_{m,m}$ generated by the $Y_{i,1}$, with $i \in \{ 1,
\dots,m \}$. 
Observe that $J=Y_{1,1}R + \dots +Y_{m,1}R$ and
$J'=Y_{1,1}R_{m,m} + \dots +Y_{m,1}R_{m,m}$. 
We have just proved that $\sigma(
Y_{j,1} ) \in J'$ for all $j \in \{ 1, \dots , m\}$. 
Hence $\sigma(J')=J'$ and
$\sigma (J)=J$. \\$ $

Thus,  
$\sigma $ induces an automorphism of $R_{m,m}/J' \simeq O_q(M_{m,m-1})$.
Now, because of the induction hypothesis (and the isomorphism $O_q(M_{m,m-1})
\simeq O_q(M_{m-1,m})$), we obtain $\sigma(Y_{i,\alpha})=\mu_{i,\alpha}
Y_{i,\alpha} + \sum_{\underline{\gamma} \in \Gamma} c_{\underline{\gamma}}
Y_{1,1}^{\gamma_{1,1}} Y_{1,2}^{\gamma_{1,2}} \dots Y_{m,m}^{\gamma_{m,m}}$
with $\mu_{i,\alpha} \neq 0$ and at least 
one of $\gamma_{k,1} > 0$ for each nonzero term of the sum. (Observe that
we can apply the inductive hypothesis since $(m-1,m) \neq (1,3)$.) \\$ $

On the other hand, it follows from Proposition \ref{imagebi2} or Proposition
\ref{imagebi2bis} that \\$\sigma(b_{n-1})= \lambda b_{n-1}$ with
$\lambda \in \mathbb{C}^*$. Moreover, recalling that $b_{n-1}=[1 \dots m
| 2 \dots m+1]$, we deduce from a
transposed version of \cite[Lemma 5.1]{glduke} that, for all $1 \leq i,\alpha
\leq m$, we have 
$$
b_{n-1} Y_{i,\alpha} = \qdot Y_{i,\alpha}
b_{n-1}
$$ 
with $\bullet = -1$ if $i=1$ and $0$ otherwise.
Let $1 \leq i,\alpha \leq m $ with $\alpha \neq 1$. Thus we must have
 $b_{n-1} \sigma(Y_{i,\alpha}) = \sigma(Y_{i,\alpha})b_{n-1} $. Since $\sigma(Y_{i,\alpha})= \mu_{i,\alpha}Y_{i,\alpha} +
\sum_{\underline{\gamma} \in \Gamma} c_{\underline{\gamma}}
Y_{1,1}^{\gamma_{1,1}} Y_{1,2}^{\gamma_{1,2}} \dots Y_{m,m}^{\gamma_{m,m}}$
with at least 
one of $\gamma_{k,1} \neq 0$ for each nonzero term of the sum, this
implies that the sum is empty and so
$\sigma(Y_{i,\alpha})=\mu_{i,\alpha}Y_{i,\alpha}$ for all $1 \leq i,\alpha
\leq m$ with $\alpha \neq 1$. \\$ $

Now $\sigma$ induces an automorphism of $R_{m,m}$ that acts on the
indeterminates $Y_{i,\alpha}$ ($1 \leq i \leq m$ and $1 < \alpha \leq m$) by
multiplication by scalars. Moreover, it follows from Proposition \ref{imagebi2}
or Proposition \ref{imagebi2bis} that $\sigma$ also acts on the 
$b_i=\left[i-n+1, \dots , m \mid 1 , \dots , m+n-i \right] $ ($n \leq i \leq
m+n-1$) by multiplication by scalars. This forces $\sigma (Y_{i,1})
= \mu_{i,1} Y_{i,1}$ for all $1 \leq i \leq m$.\\$ $

This establishes the following claim:

\begin{claim}
\label{carre}
For all $(i,\alpha) \in \gc 1, m \dc^2$, there exists $\mu_{i,\alpha} \in \mathbb{C}^*$ such that $\sigma(Y_{i,\alpha})= \mu_{i,\alpha}Y_{i,\alpha}$.
\end{claim}

It remains to consider the case that $\alpha >m$. Let us now distinguish between two cases.

If $n =m+1$, then $\sigma$ is an automorphism of $R=O_q(M_{m,m+1})$ that acts
on the indeterminates $Y_{i,\alpha}$ ($1 \leq i,\alpha \leq m$) by
multiplication by scalars. Further, because of Proposition \ref{imagebi2} or
Proposition \ref{imagebi2bis}, $\sigma$ also acts on the quantum minors
$b_i=\left[1, \dots , i \mid n-i+1 ,\dots , n \right]$ ($1 \leq i \leq m$) by
multiplication by scalars. This forces $\sigma (Y_{i,m+1}) =
\mu_{i,m+1} Y_{i,m+1}$ for all $1 \leq i \leq m$. This finishes the proof of
Theorem \ref{automxp} in the case where $n=m+1$. \\$ $

Now, assume that $m+1 < n$. Recall that  $J$ 
denotes the two-sided ideal generated by the
$Y_{k,1}$, $1 \leq k \leq m$. Note that $J=Y_{1,1}R + \dots +Y_{m,1}R$. It is
well-known that the monomials $$Y_{1,1}^{\gamma_{1,1}}Y_{2,1}^{\gamma_{2,1}}
\dots Y_{m,1}^{\gamma_{m,1}} \dots Y_{1,n}^{\gamma_{1,n}} \dots
Y_{m,n}^{\gamma_{m,n}}$$ form a PBW basis of $R$. 
Observe that
$x=\sum_{\underline{\gamma} \in \Gamma} c_{\underline{\gamma}}
Y_{1,1}^{\gamma_{1,1}} \dots Y_{m,1}^{\gamma_{m,1}} \dots
Y_{1,n}^{\gamma_{1,n}} \dots Y_{m,n}^{\gamma_{m,n}}$ belongs to $J$ if and
only if at least 
one of the $\gamma_{k,1}\geq 1$ for all $\gamma \in \Gamma$ such that 
$ c_{\underline{\gamma}} \neq 0$.

It follows from Claim \ref{carre} that $\sigma(J)=J$. Hence $\sigma$ induces
an automorphism of $R/J \simeq {\cal O}_q(M_{m,n-1})$. Since $m < n- 1$ and $m
\geq 2$, it follows from the inductive hypothesis that, for all $(i,\alpha)$
with $m+1 \leq \alpha \leq n$, we can write
\begin{eqnarray}
\label{equationfinal1} 
\sigma(Y_{i,\alpha}) & = & \mu_{i,\alpha}Y_{i,\alpha}+\sum_{\underline{\gamma}
\in \Gamma} c_{\underline{\gamma}} Y_{1,1}^{\gamma_{1,1}} \dots
Y_{m,1}^{\gamma_{m,1}} \dots Y_{1,n}^{\gamma_{1,n}} \dots
Y_{m,n}^{\gamma_{m,n}}
\end{eqnarray}
where $\Gamma$ is a finite subset of $\mathbb{N}^{mn}$ and at least 
one of the
$\gamma_{k,1} \geq 1$ for all $\gamma \in \Gamma$ such that 
$ c_{\underline{\gamma}} \neq 0$. \\$ $

Let $K$ the two-sided ideal generated by the $Y_{k,2}$, $1 \leq k \leq m$.
Note that an element $$x=\sum_{\underline{\gamma} \in \Gamma'}
c_{\underline{\gamma}} Y_{1,1}^{\gamma_{1,1}} \dots Y_{m,1}^{\gamma_{m,1}}
\dots Y_{1,n}^{\gamma_{1,n}} \dots Y_{m,n}^{\gamma_{m,n}} \in R$$ belongs to
$K$ if and only if at least 
one of the $\gamma_{k,2} \geq 1$ for all $\gamma \in
\Gamma'$ such that 
$ c_{\underline{\gamma}} \neq 0$.

It follows from Claim \ref{carre} that $\sigma(K)=K$. Hence $\sigma$ induces
an automorphism of $R/K \simeq {\cal O}_q(M_{m,n-1})$. Since $m < n- 1$ and $m
\geq 2$, the inductive hypothesis applies and, for all $(i,\alpha)$
with $m+1 \leq \alpha \leq n$, we can write
\begin{eqnarray}
\label{equationfinal2}
\sigma(Y_{i,\alpha}) & = &
\mu'_{i,\alpha}Y_{i,\alpha}+\sum_{\underline{\gamma} \in \Gamma'}
c_{\underline{\gamma}} Y_{1,1}^{\gamma_{1,1}} \dots Y_{m,1}^{\gamma_{m,1}}
\dots Y_{1,n}^{\gamma_{1,n}} \dots Y_{m,n}^{\gamma_{m,n}}
\end{eqnarray}
where $\Gamma'$ is a finite subset of $\mathbb{N}^{mn}$ and at least 
one of the
$\gamma_{k,2} \geq 1$ for all $\gamma \in \Gamma'$ such that 
$ c_{\underline{\gamma}} \neq 0$.
 \\$ $

Let $(i,\alpha)$ with $m+1\leq \alpha \leq n$. It remains to prove that
$\sigma$ acts on $Y_{i,\alpha}$ by multiplication by a scalar. First, identifying the two expressions (\ref{equationfinal1}) and
(\ref{equationfinal2}) of $ \sigma(Y_{i,\alpha})$ in the PBW basis of $R$
leads to:
\begin{eqnarray}
\label{equationfinal3} 
\sigma(Y_{i,\alpha}) 
& = & \mu_{i,\alpha}Y_{i,\alpha}+\sum_{\underline{\gamma} \in \Gamma} c_{\underline{\gamma}} Y_{1,1}^{\gamma_{1,1}} \dots  Y_{m,1}^{\gamma_{m,1}} \dots Y_{1,n}^{\gamma_{1,n}} \dots Y_{m,n}^{\gamma_{m,n}}
\end{eqnarray}
where $\Gamma$ is a finite subset of $\mathbb{N}^{mn}$ such that 
at least 
one of the
$\gamma_{k,1} \geq 1$  and at least 
one of the $\gamma_{l,2}
\geq 1$ for all $\gamma \in \Gamma$ 
such that 
$ c_{\underline{\gamma}} \neq 0$.

By Proposition \ref{imagebi2} (if $n \neq 3m$) or Proposition
\ref{imagebi2bis} (if $n = 3m$), 
the automorphism $\sigma$ acts on $b_m=[1 \dots m \mid n-m+1
\dots n ]$ by multiplication by a scalar. Further, it follows from a
transposed version of \cite[Lemma 5.1]{glduke} that
$$
b_m Y_{i,\alpha}= \left\{ \begin{array}{ll} 
q^{-1} Y_{i,\alpha} b_m & \mbox{ if } \alpha \leq n-m \\
Y_{i,\alpha} b_m & \mbox{ if } \alpha \geq n-m+1.
\end{array}
\right.
$$
Hence
$$
b_m \sigma( Y_{i,\alpha} ) = \left\{
\begin{array}{ll}
q^{-1} \sigma(Y_{i,\alpha}) b_m & \mbox{ if }
\alpha \leq n-m \\
\sigma(Y_{i,\alpha}) b_m & \mbox{ if } \alpha \geq
n-m+1,
\end{array}
\right.
$$
Thus
$$
\sum_{\underline{\gamma} \in \Gamma} c_{\underline{\gamma}}
q^{-\gamma_{1,1}-\dots-\gamma_{m,1}-\dots -\gamma_{1,n-m}- \dots
-\gamma_{m,n-m}} Y_{1,1}^{\gamma_{1,1}} \dots Y_{m,1}^{\gamma_{m,1}} \dots
Y_{1,n}^{\gamma_{1,n}} \dots Y_{m,n}^{\gamma_{m,n}}
$$
$$
= \qdot \sum_{\underline{\gamma} \in \Gamma} c_{\underline{\gamma}}
Y_{1,1}^{\gamma_{1,1}} \dots Y_{m,1}^{\gamma_{m,1}} \dots
Y_{1,n}^{\gamma_{1,n}} \dots Y_{m,n}^{\gamma_{m,n}},
$$
where $\qdot = -1$ if $ \alpha \leq n-m$, and $\qdot=0$ otherwise.

Consequently, 
\begin{eqnarray}
\label{fin2}
\gamma_{1,1}+\dots+\gamma_{m,1}+\dots +\gamma_{1,n-m}+ \dots +\gamma_{m,n-m} &
= & \left\{ \begin{array}{ll} 1 & \mbox{ if } \alpha \leq n-m \\ 0 & \mbox{ if
} \alpha \geq n-m+1,
\end{array}
\right.
\end{eqnarray} 
for all $\underline{\gamma} \in
\Gamma$ such that $c_{\underline{\gamma}}\neq 0$, since $q$ is not a root of
unity.
$ $

On the other hand, recall that
$$\sigma(Y_{i,\alpha})=\mu_{i,\alpha}Y_{i,\alpha}+\sum_{\underline{\gamma} \in
\Gamma} c_{\underline{\gamma}} Y_{1,1}^{\gamma_{1,1}} \dots
Y_{m,1}^{\gamma_{m,1}} \dots Y_{1,n}^{\gamma_{1,n}} \dots
Y_{m,n}^{\gamma_{m,n}}$$
where, for each nonzero term of the sum, at least 
one of the $\gamma_{k,1}$ is a
positive integer and at least 
one of the $\gamma_{l,2}$ is a positive integer. 
Hence,
$\gamma_{1,1}+\dots+\gamma_{m,1}+\dots +\gamma_{1,n-m}+
\dots +\gamma_{m,n-m} \geq 2$ 
for all $\underline{\gamma} \in \Gamma$ with $c_{\underline{\gamma}}\neq 0$, 
since $n-m \geq 2$. 
This
contradicts (\ref{fin2}). Thus $\Gamma$ must be empty for all $(i,\alpha)$
with $\alpha \geq  m$ 
and so $\sigma(Y_{i,\alpha}) = \mu_{i,\alpha} Y_{i,\alpha} $
for all $\alpha > m$. This finishes the proof.
\end{proof}

\begin{corollary}
Assume that $m < n$ and $(m,n) \neq (1,3)$. Let $\sigma$ be an automorphism of $\oqmmn$. Then there exist unique nonzero complex numbers 
$h_1,\dots,h_m,h'_1,\dots,h'_{n-1}$ such that $\sigma(Y_{i,\alpha})= h_ih'_{\alpha} Y_{i,\alpha}$ (with the convention $h'_n=1$).
\\Thus, $\aut (\oqmmn)$ is isomorphic to the torus $(\mathbb{C}^*)^{m+n-1}$.
\end{corollary}
\begin{proof}
Let $\sigma$ be an automorphism of $\oqmmn$. By Theorem \ref{automxp},
there exists a family $(\mu_{i,\alpha})_{(i,\alpha) \in \gc 1,m \dc \times \gc
1,n \dc}$ of elements of $\mathbb{C}^*$ such that $\sigma(Y_{i,\alpha})=
\mu_{i,\alpha}Y_{i,\alpha}$ for all $(i,\alpha) \in \gc 1, m \dc \times \gc
1,n \dc$. Recall that, if $(i,\alpha),(j,\beta) \in \gc 1, m \dc \times \gc
1,n \dc$ with $i <j$ and $\alpha <\beta$, then $Y_{j,\beta}Y_{i,\alpha}=
Y_{i,\alpha} Y_{j,\beta} -(q-q^{-1}) Y_{i,\beta} Y_{j,\alpha}$. Hence, since
$\sigma$ is an automorphism, we must have $\mu_{i,\alpha}\mu_{j,\beta} =
\mu_{i,\beta} \mu_{j,\alpha}$ for all $(i,\alpha),(j,\beta) \in \gc 1, m \dc
\times \gc 1,n \dc$ with $i <j$ and $\alpha <\beta$. In other words, the
matrix $(\mu_{i,\alpha})_{(i,\alpha) \in \gc 1,m \dc \times \gc 1,n \dc}$ has
rank $1$. Hence, there exist unique nonzero complex numbers
$h_1,\dots,h_m,h'_1,\dots,h'_{n-1},h'_n=1$ such that $\mu_{i,\alpha}= h_i h'_\alpha$
for all $i,\alpha$, as desired.
\end{proof}

Note that, in the exceptional case where $(m,n) = (1,3)$, 
the automorphism group of
$\oqmm13$ has been computed by Alev and Chamarie, \cite[Th\'eor\`eme
1.4.6]{alevchamarie}; in this case, the group is not isomorphic to the torus
$(\mathbb{C}^*)^3$, since the second case of Lemma~\ref{imagebi1bis} does
arise. Indeed, Alev and Chamarie show that any automorphism of $\oqmm13$ is of
the form 
$$ 
\sigma(Y_{11}) = \mu_1 Y_{11}, \quad \sigma(Y_{12}) = \mu_2 Y_{12}
+ \lambda Y_{11}Y_{13}, \quad \sigma(Y_{13}) = \mu_3 Y_{13},  
$$ 
where $\mu_i \neq 0$ and $\lambda$ are complex numbers; 
so that the
automorphism group is isomorphic to the semidirect 
product $\mc \times (\mc^*)^3$. This result can easily be obtained from 
our analysis. \\

\noindent {\bf Acknowledgments.} We thank Jacques Alev, Ken Goodearl and Laurent
Rigal for helpful conversations and comments.




\newpage

\noindent S Launois:\\
School of Mathematics, University of Edinburgh,\\
James Clerk Maxwell Building, King's Buildings, Mayfield Road,\\
Edinburgh EH9 3JZ, Scotland\\
E-mail: stephane.launois@ed.ac.uk \\

\noindent T H Lenagan: \\
School of Mathematics, University of Edinburgh,\\
James Clerk Maxwell Building, King's Buildings, Mayfield Road,\\
Edinburgh EH9 3JZ, Scotland\\
E-mail: tom@maths.ed.ac.uk \\


\end{document}